\documentclass[11pt,a4paper]{article}
\usepackage{a4wide}
\usepackage{theorem}
\usepackage{amsmath,amssymb,calc}
\usepackage[english]{babel}
\usepackage{graphicx,array}
\usepackage[numbers]{natbib}
\usepackage{tikz}
\usepackage{color} 
\usepackage{enumerate}
\usepackage[toc,page]{appendix}
\usepackage[charter]{mathdesign}
\usepackage{algorithm}
\usepackage{algpseudocode}
\usepackage{hyperref}
\hypersetup{
    bookmarksopen=false,
    bookmarksnumbered=true,
    pdftitle={Networks of fixed-cycle intersections},
    pdfauthor={M.A.A. Boon, J.S.H. van Leeuwaarden}
}
\ifpdf
  \hypersetup{colorlinks=true,linkcolor=black,urlcolor=black,citecolor=black,pdfpagemode=UseOutlines,plainpages=false,pdfpagelabels}
\else
  \hypersetup{colorlinks=false}
\fi

\newtheorem{theorem}{Theorem}[section]
\newtheorem{assumption}[theorem]{Assumption}

\newtheorem{lemma}[theorem]{Lemma}
\newtheorem{remark}[theorem]{Remark}

\numberwithin{equation}{section}

\newenvironment{rem}{\begin{remark}\normalfont}{\end{remark}}
\newenvironment{ass}{\begin{assumption}\normalfont}{\end{assumption}}

\newcommand{\ee}{\text{e}}
\newcommand{\E}{\mathbb{E}}
\renewcommand{\P}{\mathbb{P}}

\newcommand{\G}{\mathcal{G}}

\newcommand{\Sk}{\mathcal{S}^{(k)}}
\newcommand{\Gk}{\mathcal{G}}
\newcommand{\equaldist}{\,{\buildrel d \over =}\,}

\newcommand{\un}{\underline{n}}
\providecommand{\href}[2]{#2}

\usetikzlibrary{decorations.pathreplacing}

\newcommand{\trafficsignal}[2][0]{
\begin{scope}[shift={#2}, rotate=#1]
\draw[black,thick,fill=black](0,0) rectangle (1,2);
\draw[fill=red](0.5,1.6) circle [radius=0.3];
\draw[fill=yellow](0.5,1.0) circle [radius=0.3];
\draw[fill=green](0.5,0.4) circle [radius=0.3];
\end{scope}
}

\newcommand{\car}[2][0]{
\begin{scope}[shift={#2}, rotate=#1, scale=0.028]
\draw[white,fill=white](0.8,14.2) -- (1.4,15.6) -- (2.2,16.6) -- (2.6,23.8) -- (3.2,28.4) -- (5.4,36.2) -- (9,42.4) -- (19.3,43) -- (30.1,43.3) -- (34.8,42.9) -- (38.9,42.1) -- (41.5,41.1) -- (45.8,38.7) -- (56,31.3) -- (65,23.6) -- (66.2,22.1) -- (67.1,19.7) -- (67.3,17.5) -- (68.3,16.4) -- (68.3,16.4) -- (69.2,14.9) -- (68.7,14) -- (60.6,7.6) -- (50.9,5.6) -- (19.5,5.6) -- (9.7,7.4) -- (1.6,13.5);
\draw[gray,fill=gray](9.6,7.5) circle [radius=10];
\draw[gray,fill=white](9.6,7.5) circle [radius=8.5];
\draw[gray,fill=gray](60.6,7.5) circle [radius=10];
\draw[gray,fill=white](60.6,7.5) circle [radius=8.5];
\draw[gray,fill=gray](48.8,27.2) -- (20,31.1) -- (17.7,32.9) -- (16.9,36) -- (18.5,39.5) -- (20.6,41.1) -- (22.4,41.6) -- (33.6,41.1) -- (37.6,40.3) -- (40.6,39) -- (49.9,33) -- (51.1,31.9) -- (51.5,30.3) -- (51,28.8) -- (50,28.1);
\end{scope}
}

\newcommand\blfootnote[1]{%
  \begingroup
  \renewcommand\thefootnote{}\footnote{#1}%
  \addtocounter{footnote}{-1}%
  \endgroup
}

\title{Networks of fixed-cycle intersections\blfootnote{Department of Mathematics and Computer Science, Eindhoven University of Technology, P.O. Box 513, 5600MB Eindhoven, The Netherlands. Email: \href{mailto:m.a.a.boon@tue.nl}{m.a.a.boon@tue.nl} and \href{mailto:j.s.h.v.leeuwaarden@tue.nl}{j.s.h.v.leeuwaarden@tue.nl}}}
\author{M.A.A. Boon
 \and J.S.H. van Leeuwaarden
 }

\begin{document}
\maketitle

\begin{abstract}
We present an algorithmic method for analyzing networks of intersections with static signaling, with as primary example a line network
that allows traffic flow over several intersections in one main direction.
The method decomposes the network into separate intersections and treats each intersection in isolation using an extension of the fixed-cycle traffic-light (FCTL) queue. The network effects are modeled by matching the output process of one intersection with the input process of the next (downstream) intersection. This network analysis provides insight into wave phenomena due to vehicles experiencing progressive cascades of green lights and sheds light on platoon forming in case of imperfections. Our algorithm is shown to match results from extensive discrete-event simulations and can also be applied to more complex network structures.
\end{abstract}

\section{Introduction}\label{sect:introduction}
Intersections are natural bottlenecks and crucially influence the dynamics of urban traffic. Traffic lights trigger a switching process meant to manage conflicting traffic flows. The coordination is sometimes done dynamically, according to sensor data of currently existing traffic flows; otherwise it is done statically, by the use of timers.
While intersections can be studied in isolation \cite{darroch,tachet2016revisiting,fctlsolo}, the larger picture
of networks of multiple intersections is increasingly important, also in view of the rapid growth of urbanization \cite{bettencourt2013origins,urban}. This paper contributes to the  theoretical underpinning of traffic networks by extending classical models for isolated intersections to models for networks of intersections with static signaling.

Think of a series of traffic lights designed to let traffic flow over several intersections in one main direction. Any vehicle traveling along (at an approximate prescribed speed) wants  to meet a progressive cascade of green lights, and not have to stop at intersections. In practical use, only a group of vehicles -- referred to as {\it platoon} -- can pass the intersection before the time band is interrupted to give way to other traffic flows. The platoon sizes are governed by the signal times. Our method to model such situations consists of two ingredients: An extension of a classical queueing model for one isolated intersection that can deal with correlated input and that allows for a detailed characterization of the output process of an intersection, and an algorithm for network analysis that decomposes the series of queues into multiple isolated queues. While interesting in their own right, and likely to find more applications in transportation science, network analysis of intersections requires the delicate combination of both ingredients. We now discuss each of them separately.


\vspace{.3cm}
\noindent{\bf Queueing model for one intersection.} The classical model for an isolated intersection that we adopt and extend in this paper is the fixed-cycle traffic-light (FCTL) queue; one of the most well-studied stochastic models in traffic engineering \cite{darroch,mcneill,newell65,webster}. Vehicles arrive to an intersection controlled by a traffic light and form a queue. The time scale is divided into time intervals of unit length, and the traffic light  alternates between red and green periods of fixed durations $r$ and $g$ time units. Delayed vehicles depart during the green period, where it takes one time unit for each delayed vehicle to depart; departures thus occur at equally spaced times until either the queue dissipates or the green phase terminates. Darroch \cite{darroch} obtained the probability generating function (pgf) of the steady-state overflow queue (the number of vehicles waiting in front of the traffic light at the end of a green period) and the pgf of the steady-state delay was obtained in van Leeuwaarden \cite{fctlsolo}. Hence, all information about the distribution of the steady-state overflow queue and steady-state delay in the FCTL queue can be obtained from the results in \cite{darroch,fctlsolo}, including all moments of the steady-state queue length and delay, and the distribution of the output process (the way vehicles leave the intersection).
The output process of the first intersection is of crucial importance for the present paper, because it will serve as input process for some other signalized intersection. Moreover, the output process of a second intersection serves as input for a third intersection, and so forth. This network effect acts as a filter that modifies, and perhaps streamlines, the arrival process at consecutive intersections. Therefore, we shall address in this paper the technical challenge of extending the classical FCTL queue to allow for nonuniform and hence time-dependent correlated arrival processes. We call this extended model the generalized FCTL queue.

\vspace{.3cm}
\noindent{\bf Network algorithm.} A network of intersections with correlated input and output processes appears not solvable. We therefore develop an approximation scheme to evaluate the system performance based on decomposition. While this approach has been successfully applied to classic queueing networks \cite{bitrandasu,kuehn,whitt1983queueing}, a network of generalized FCTL queues poses additional challenges due to the non-synchronized cyclic structures and inherently correlated arrival processes.
We decompose the network into isolated generalized FCTL queues, which are then analyzed separately by assuming specific arrival processes, and in particular the output process of one intersection serves as the input process of an upstream intersection, hence creating the correlation structure that comes with network topologies.

\vspace{.3cm}
\noindent{\bf Outline of the paper.}
In Section~\ref{sec:queuelengths} we provide a detailed model description of the generalized FCTL queue. In Section~\ref{sec:main} we present the full analytic solution of the generalized FCTL queue, both in terms of a formal characterization of the probability generating functions of the queue length distribution, and in terms of practically implementable algorithms for calculating the queue length distribution, for any given correlated arrival pattern.
In Section~\ref{sec:networks} we design the network algorithm based on decomposition and the results in Section \ref{sec:main}. We also compare our analytical results with extensive discrete-event simulation of the same network model. In Section \ref{sec:conclusions} we present conclusions.

\section{FCTL queue with correlated arrivals}\label{sec:queuelengths}

We now present a generalization of the classical FCTL queue that can deal with correlated arrivals.
In Subsection \ref{s31} we detail the model and its underlying assumptions. We then discuss in Subsection \ref{s32} an example of an arrival pattern that contains some crucial features than are anticipated in network settings. The numerical calculations for that example were performed with the algorithmic method developed in Section~\ref{sec:main}.

\subsection{Model description}\label{s31}

The first two model assumptions are adopted from the classical FCTL queue \cite{fctlsolo}:

\begin{ass}\label{discrete-time} (Discrete-time assumption) The time axis is divided into constant time intervals of unit length, so-called slots, where each slot corresponds to the time needed for a delayed vehicle to depart from the queue. The green and red periods, and thus the cycle time $c$, are assumed to be fixed multiples of one slot. Hence, $g,r,c$ are integers expressed in slots. Those vehicles that arrive to the queue and are delayed, join the queue at the end of the slot in which they arrive.
\end{ass}

\begin{ass}\label{fctlassumption} (FCTL assumption) For those cycles in which the queue clears before the green period terminates, all vehicles that arrive during the residual green period pass through the system and experience no delay whatsoever.
\end{ass}

The FCTL assumption lets vehicles that arrive during the residual green period pass the intersection without slowing down, and therefore the discharge rate of these vehicles is larger than the discharge rate of the delayed vehicles (one per time unit). Because of the huge difference in discharge rates of delayed vehicles (these vehicles have to accelerate) and non-delayed vehicles, the FCTL assumption is a sensible assumption. The next assumption is new:

\begin{ass}\label{ia} (Correlated arrivals assumption)  Let $Y_{i,n}$ denote the number of vehicles that arrive to the intersection during slot $i$ in cycle $n$. The random variables $Y_{i,n}$ are allowed to be dependent within cycle $n$, but we assume that $Y_{i,n}$ and $Y_{i,m}$ are independent when $n\neq m$.
\end{ass}

Notice that Assumptions \ref{discrete-time} and \ref{ia} together make that the queue lengths at the end of time slots can be modeled as a discrete-time Markov chain. This feature is exploited in Section \ref{sec:main} to find a fully analytic characterization for the steady-state queue-length distribution.

It is crucial that Assumption \ref{ia} is less restrictive than its counterpart in \cite{fctlsolo} that assumes the $Y_{i,n}$ to be independent and identically distributed (i.i.d.).~We need to move beyond this i.i.d.~assumption in order to consider the correlated $Y_{i,n}$ sequences as they can occur in real network settings, for instance when  the output process of one intersection (or FCTL queue) forms the input process for another intersection. Because the first intersection alters the original arrival pattern, the second intersection is likely to be confronted with
platoons of vehicles that have been delayed by the upstream red signal.

The generalized FCTL queue defined by Assumptions \ref{discrete-time}-\ref{ia} is in essence a queueing system with multiple customer types and batch arrivals, where $Y_{i,n}$ can be interpreted as the number of type $i$ customers in batch $n$. Denote the pgf of the joint distribution of $(Y_{1,n}, Y_{2,n}, \dots, Y_{c,n})$ by
\[
Y_n(y_{1}, y_{2}, \dots, y_{c}) = \E\Big[\prod_{i=1}^c y_i^{Y_{i,n}}\Big].
\]
Although arrivals within a batch can be correlated, successive batches are i.i.d. Multi-type queueing models with batch arrivals have been well-studied \cite{henderson1990,laevens,pollingbatch}. A marked difference however, is that arrivals in the FCTL queue do not join the queue instantaneously, but are dictated when to arrive according to their type: type $i$ arrivals join the queue in time slot $i$, for $i=1,2,\dots,c$. Another crucial difference is Assumption \ref{fctlassumption}, of course, which is very specific for traffic-light settings.

\subsection{Motivating example}\label{s32}
We now give an exemplary arrival pattern with features that are anticipated in network settings.
We consider a generalized fixed-cycle traffic light queue with traffic arriving from one synchronized upstream traffic intersection.
The cyclic arrival pattern at this queue, illustrated in Figure \ref{fig:isolatedintersection}, starts with the arrival of a platoon of delayed vehicles from the major upstream flow. The settings are synchronized such that the signal turns green at the exact moment that the first car in this platoon arrives.
This platoon is followed by a phase of free flow, with arriving     vehicles that were not delayed at the upstream intersection. 
After a short period of no arrivals, we have a similar pattern of a platoon followed by free flow, arriving from a minor upstream flow.
A more detailed description of the arrival process can be found in Appendix \ref{appendix:input}.

In Section \ref{sect:networkexamples} we show that this specific arrival pattern arises naturally in network settings. 
Arrivals in the same cycle are correlated, but we assume independence between arrivals in successive cycles. 

\begin{figure}[ht]
\begin{center}
\begin{tikzpicture}[scale=0.3]
\draw[black,fill=red](0,0) rectangle (40,2);
\draw[black,fill=green](0,0) rectangle (20,2);
\draw[black](2,0) rectangle (2,2);
\draw[black](4,0) rectangle (4,2);
\draw[black](6,0) rectangle (6,2);
\draw[black](8,0) rectangle (8,2);
\draw[black](10,0) rectangle (10,2);
\draw[black](12,0) rectangle (12,2);
\draw[black](14,0) rectangle (14,2);
\draw[black](16,0) rectangle (16,2);
\draw[black](18,0) rectangle (18,2);
\draw[black](22,0) rectangle (22,2);
\draw[black](24,0) rectangle (24,2);
\draw[black](26,0) rectangle (26,2);
\draw[black](28,0) rectangle (28,2);
\draw[black](30,0) rectangle (30,2);
\draw[black](32,0) rectangle (32,2);
\draw[black](34,0) rectangle (34,2);
\draw[black](36,0) rectangle (36,2);
\draw[black](38,0) rectangle (38,2);
\node[below] at (5,0) {platoon 1};
\node[below] at (17,0) {free flow};
\node[below] at (31,0) {platoon 2};
\node[below] at (36.8,0) {free flow};
\draw[decorate, decoration={brace}, yshift=2ex]  (0,2) -- node[above=0.4ex] {$B_1$}  (13.9,2);
\draw[decorate, decoration={brace}, yshift=2ex]  (14.1,2) -- node[above=0.4ex] {$10-B_1$}  (20,2);
\draw[decorate, decoration={brace}, yshift=2ex]  (30,2) -- node[above=0.4ex] {$B_2$}  (33.9,2);
\draw[decorate, decoration={brace}, yshift=2ex]  (34.1,2) -- node[above=0.4ex] {\ \ $3-B_2$}  (36,2);
\foreach \x in {1,...,20}
    \node[below] at (2*\x-1,1.7) {\x};
\end{tikzpicture}

\end{center}
\caption{The arrival pattern in Section \ref{s32}. 
The lengths of platoons 1 and 2 are respectively $B_1$ and $B_2$. The cycle length is $c=g+r=10+10=20$ time slots.
}
\label{fig:isolatedintersection}
\end{figure}

\begin{figure}[ht]
\begin{center}
\includegraphics[width=0.5\textwidth]{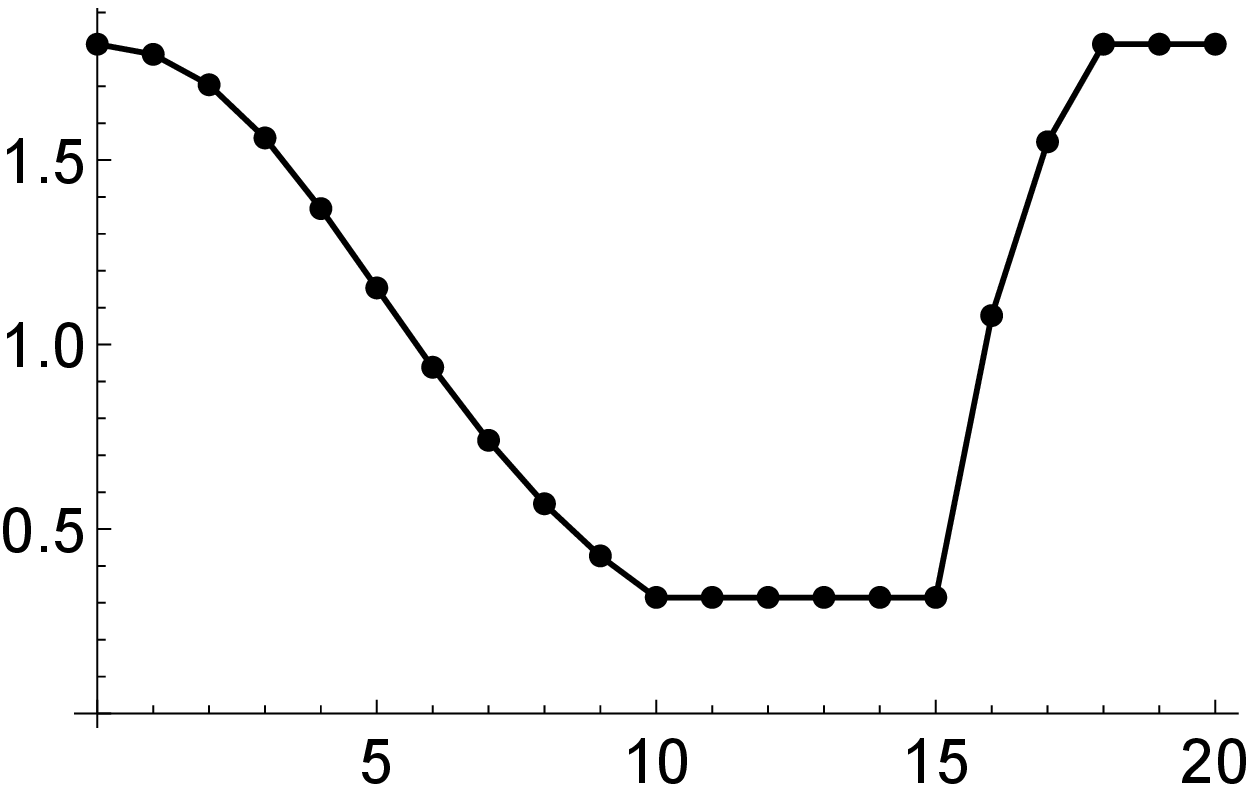}
\begin{picture}(0,0)
\put(0,0.3){ $i$}
\put(-8,5){ $\E[X_i]$}
\end{picture}
\end{center}
\caption{The queue length of the isolated intersection in Section \ref{s32}. The black dots represent queue lengths obtained via a microscopic simulation of our model and the solid black line  is obtained by the algorithm developed in Section \ref{sec:main}. }
\label{fig:isolatedintersectionQL}
\end{figure}
Figure \ref{fig:isolatedintersectionQL} shows the mean queue length during a cycle. The results from a microscopic discrete-event simulation of this generalized FCTL queue have been used to validate the analytic method. The longest queue occurs just before the traffic light turns green, and the queue is typically minimized at the end of green periods. The graph shows an interesting pattern, where the queue length remains constant for the first part of the red period, due to the absence of arrivals during that period. The arrival of platoon~2 results in a steep increase in queue length, starting after time slot~15. 
The shape of the graph in Figure~\ref{fig:isolatedintersectionQL} is different from the standard FCTL queue, which is typically more V-shaped due to the constant arrival rate throughout the cycle.

The analytic results obtained in the next section allow us to determine the queue-length distribution. Table~\ref{tbl:queuelengthprobsExample1} shows probabilities of the queue length exceeding certain levels during the cycle. We see that the probability that the queue exceeds five vehicles at the cycle start is $0.015$. At the end of the green period this probability is only $0.002$. At any arbitrary moment, the probability that there are more than five delayed vehicles is $0.008$.

\begin{table}[!ht]
\[
\begin{array}{l|cccccc}
\hline
X & \P(X\geq1) & \P(X\geq2) & \P(X\geq3) & \P(X\geq4) & \P(X\geq5) & \P(X\geq6)  \\
\hline
X_0     & 0.829 & 0.547 & 0.302 & 0.075 & 0.036 & 0.015  \\
X_{10}  & 0.159 & 0.089 & 0.042 & 0.014 & 0.006 & 0.002  \\
\bar{X} & 0.496 & 0.294 & 0.146 & 0.042 & 0.019 & 0.008  \\
\hline
\end{array}
\]
\caption{The queue length distribution at the beginning of a cycle $(X_0)$, at the end of a green period $(X_{10})$, and at arbitrary moments $(\bar{X})$. The notation in brackets will be introduced in Section \ref{sec:main}.}
\label{tbl:queuelengthprobsExample1}
\end{table}

\section{Analytic solution and algorithms}\label{sec:main}
We shall now show that the  generalized FCTL queue with Assumptions \ref{discrete-time}-\ref{ia} is analytically solvable. Using generating function techniques and complex analysis, we are able to obtain an explicit, analytic characterization for the steady-state queue-length distribution. In order to do so, we exploit the Markovian nature of the queueing model at the level of cycles, and account for the correlations that occur within the cycle. The recursion relation that connects consecutive queue lengths during a cycle is given in Subsection \ref{s41}, along with some further preliminaries. Then in Subsection \ref{s42} we derive the generating function for the steady-state queue lengths at the beginning of time slots throughout a cycle. In particular, we leverage this explicit  generating function to create an algorithm for calculating the complete queue-length distribution, at all points in time. Special attention is paid to the mean queue lengths in Subsection \ref{s43}.

\subsection{Preliminaries}\label{s41}
Let $X_{i,n}$ denote the queue length at the end of the $i$-th time slot in the $n$-th cycle with $n=1, 2, \dots$. Let $Y_{i,n}$ denote the number of arrivals in time slot $i$ during the $n$-th cycle, as defined in Assumption \ref{discrete-time}.

For convenience, define $X_{0,n+1}(z) := X_{c,n}(z)$ as the pgf of the queue-length distribution at the start of cycle $n+1$. We will determine the the steady-state queue length, in the limit as $n\rightarrow\infty$, at the end of each time slot. Without loss of generality, we assume that the first $g$ time slots in the cycle have a green signal, and the last $r$ time slots are red. The queue length in this slot-based system then evolves according to
\begin{equation}
X_{k+1,n}=\begin{cases}
X_{k,n}+Y_{k+1,n}-1 & \qquad \text{if } X_{k,n}>0 \text{ and }k = 0, \dots, g-1,\\
0 & \qquad \text{if } X_{k,n}=0 \text{ and }k = 0, \dots, g-1,\\
X_{k,n}+Y_{k+1,n} & \qquad \text{for }k = g, \dots, c-1.
\end{cases}
\label{eqn:queuelengthevolutionold}
\end{equation}

In the classical FCTL queue \cite{darroch,fctlsolo}, with i.i.d.~$Y_{1,n},\dots,Y_{c,n}$, the function $X_{k+1,n}(z)$ can be expressed in terms of $X_{k,n}(z)$ by conditioning on the event that $X_{k,n}$ is equal to zero, or not. For the generalized FCTL queue more detailed information is needed, which is why we will use a lattice path counting approach to describe the queue-length evolution throughout the cycles.

Let $\mathbf{1}_{\{A\}}$ be the indicator function for event $A$. Define the function
\[
T(x_0,n_1,n_2,\dots,n_k) = \mathbf{1}_{\{\min(x_0, x_0+n_1-1, x_0+n_1+n_2-2, \dots, x_0+\sum_{i=1}^k n_i-k)\leq 0\}},
\]
for $k=0, 1, 2, \dots, g$. The parameters $x_0, n_1, \dots,n_k$ are allowed to be random variables. We define $\Sk \subset \mathbb{N}^k$ as the set of $k$-dimensional integer-valued vectors $(l, n_1, n_2, \dots, n_{k-1})$ for which $T(l, n_1, n_2, \dots, n_{k-1})=1$.
For $j=0,1,\dots, k-1$, let $\Sk_j \subseteq \Sk$ denote the subset of $\Sk$ with elements $(l, n_1, n_2, \dots, n_{k-1})$ that satisfy
for $m=0, 1, \dots, j-1$,
\begin{align*}
l+\sum_{i=1}^{m} n_i-m > 0, &\qquad  \text{ and }\\
l+\sum_{i=1}^{j} n_i-j = 0.&
\end{align*}
Notice that $\Sk$ contains all possible combinations $(X_0, Y_1, Y_2, \dots, Y_{k-1})$ that will cause $X_{k-1}$ to be zero. Out of all these combinations, $\Sk_j$ contains the combinations for which $X_{j}=0$, while $X_0, \dots, X_{j-1}$ all take strictly positive values. Adopting the terminology from \cite{broek}, we say that $\Sk_j$ contains all elements for which the \emph{effective green time} equals $j$, which means to say that the queue becomes empty after precisely $j$ time slots. We will subdivide these subsets further by conditioning on their first element, which represents the queue length at the beginning of the cycle. Define $\Sk_{j,l} \subseteq \Sk_j$ as the elements of $\Sk_j$ that have first entry $l$. Note that $\Sk_{j,l}$ only needs to be defined for $l=0, 1, \dots, k-1$, for otherwise the queue cannot be empty at the end of time slot $k-1$.
Finally, define 
\[
\G_{j,l}:=\{(n_1,\dots,n_j)\, | \,(l,n_1, n_2, \dots, n_j, n_{j+1},\dots,n_k) \in \mathcal{S}^{(j+1)}_{j,l}\}
\]and notice that
\[(n_1, n_2, \dots, n_j) \in \Gk_{j,l}\]
implies
\[
(l,n_1, n_2, \dots, n_j, n_{j+1},\dots,n_k) \in \Sk_{j,l},\]
for any $n_{j+1},\dots,n_k \in \mathbb{N}$ and $k=j+1,\dots,g-1$.

\begin{rem}
A technical issue is that we allow $j$ to be zero in the definition of $\G_{j,l}$. The only way to obtain an effective green time of length zero is when the queue length was zero at the cycle start. As a consequence, we have
\[\G_{0,0}=\{(\,)\}, \qquad \G_{0,l} = \emptyset \text{ for }l=1,2,\dots ,\]
where we have used $(\,)$ to denote a vector of dimension 0.
\end{rem}

We can now express the queue lengths $X_{k,n}$ in terms of $X_{0,n}$. For $k=1,\dots,g$,
\begin{equation}
X_{k,n}=\begin{cases}
0 & \qquad\text { if }T(X_{0,n}, Y_{1,n}, Y_{2,n}, \dots, Y_{k-1,n})=1,\\
X_{0,n}+\sum_{i=1}^k Y_{k,n}-k & \qquad\text{ otherwise.}
\end{cases}
\label{eqn:queuelengthevolution1}
\end{equation}
and for $k=g+1,\dots,c$,
\begin{equation}
X_{k,n}=\begin{cases}
\sum_{i=g+1}^k Y_{k,n} & \qquad\text { if }T(X_{0,n}, Y_{1,n}, Y_{2,n}, \dots, Y_{g-1,n})=1,\\
X_{0,n}+\sum_{i=1}^k Y_{k,n}-g & \qquad\text{ otherwise.}
\end{cases}
\label{eqn:queuelengthevolution2}
\end{equation}

\subsection{Generating function solution}\label{s42}
From the recursion relation \eqref{eqn:queuelengthevolution1}, and using the fact that $X_0$ is independent of the future arrivals, we find
\begin{align*}
X_k(z)=&\sum_{l=0}^\infty \sum_{n_1=0}^\infty\dots\sum_{n_k=0}^\infty \big(1-T(l,n_1,\dots,n_{k-1})\big) z^{l+n_1+n_2+\dots+n_k-k}
\P(X_0=l,Y_1=n_1,\dots,Y_k=n_k)\\
&+\sum_{l=0}^\infty \sum_{n_1=0}^\infty\dots\sum_{n_k=0}^\infty T(l,n_1,\dots,n_{k-1})
\P(X_0=l,Y_1=n_1,\dots,Y_k=n_k)\\
=&\sum_{l=0}^\infty \sum_{n_1=0}^\infty\dots\sum_{n_k=0}^\infty z^{l+n_1+n_2+\dots+n_k-k}
\P(X_0=l,Y_1=n_1,\dots,Y_k=n_k)\\
&+\sum_{l=0}^\infty \sum_{n_1=0}^\infty\dots\sum_{n_k=0}^\infty T(l,n_1,\dots,n_{k-1})\left(1-z^{l+n_1+n_2+\dots+n_k-k}\right)
\P(X_0=l,Y_1=n_1,\dots,Y_k=n_k)
\end{align*}
Using the definitions of $\Sk_j$ and the independence between $X_0$ and $Y_1, \dots, Y_k$,
\begin{align*}
&X_k(z)= X_0(z)\frac{Y(\overbrace{z, \dots, z}^{k}, 1, \dots, 1)}{z^k}\\
&+\sum_{l=0}^{k-1}\sum_{j=0}^{k-1} \sum_{(l,n_1,\dots,n_{k-1})\in \Sk_j}\sum_{n_k=0}^\infty\left(1-z^{l+n_1+n_2+\dots+n_k-k}\right)\P(X_0=l,Y_1=n_1,\dots,Y_k=n_k).
\end{align*}
Since $\P(X_0=l,Y_1=n_1,\dots,Y_k=n_k)=q_l\P(Y_1=n_1,\dots,Y_k=n_k)$ we get
\begin{align*}
&X_k(z)
= X_0(z)\frac{Y(\overbrace{z, \dots, z}^{k}, 1, \dots, 1)}{z^k}\\
&+\sum_{l=0}^{k-1}q_l\sum_{j=0}^{k-1}
\sum_{(n_1,\dots,n_{j})\in \Gk_{j,l}}\sum_{n_{j+1}=0}^\infty\dots\sum_{n_k=0}^\infty\left(1-z^{l+n_1+n_2+\dots+n_k-k}\right)\P(Y_1=n_1,\dots,Y_k=n_k).
\end{align*}
Now, for an arbitrary element $(n_1,\dots,n_{j})$ in $\Gk_{j,l}$, it should hold that $l+n_1+\dots+n_j=j$. We can use this to write
\begin{align*}
\sum_{n_{j+1}=0}^\infty&\dots\sum_{n_k=0}^\infty\left(1-z^{l+n_1+n_2+\dots+n_k-k}\right)\P(Y_1=n_1,\dots,Y_k=n_k)\\
&=\P(Y_1=n_1,\dots,Y_{j}=n_j)-\sum_{n_{j+1}=0}^\infty\dots\sum_{n_k=0}^\infty z^{j+n_{j+1}+\dots+n_k-k}\P(Y_1=n_1,\dots,Y_k=n_k)\\
&=:f_k(1, n_1, \dots, n_j)-f_k(z, n_1, \dots, n_j),
\end{align*}
where
\begin{equation}
f_k(z, n_1, \dots, n_j)=\frac{\partial}{\partial y_1^{n_1}\partial y_2^{n_2}\dots \partial y_j^{n_j}} \,
 \frac{Y(y_1,y_2,\dots,y_j,\overbrace{z,\dots,z}^{k-j},1,\dots,1)}{n_1!n_2!\dots n_j!\,z^{k-j}}\Big|_{y_1=\dots=y_j=0}.
\label{eqn:fk}
\end{equation}
To see the latter, note that
by definition
\[
f_k(z, n_1, \dots, n_j)=\sum_{n_{j+1}=0}^\infty \sum_{n_{j+2}=0}^\infty \dots \sum_{n_{k}=0}^\infty z^{n_{j+1}+n_{j+2}+\dots +n_{k}-(k-j)}\P\left(Y_{1}=n_1, Y_{2}=n_2,\dots,Y_k=n_{k}\right),
\]
and hence
\[
f_k(1, n_1, \dots, n_j)=\P\left(Y_{1}=n_1, Y_{2}=n_2,\dots,Y_j=n_{j}\right).
\]

We thus obtain the following lemma, where we have omitted the subscript $n$ that refers to a particular cycle.
\begin{lemma} \label{thm:xkx0}
Let $\un = (n_1,\dots,n_j)$ and $q_l = \P(X_0=l)$ for $l=0, 1, \dots, g-1$. Then, for $k=1, 2, \dots, g$,
\begin{equation}
X_k(z) = X_0(z)\frac{Y(\overbrace{z, \dots, z}^{k}, 1, \dots, 1)}{z^k} +
\sum_{l=0}^{k-1}q_l \sum_{j=0}^{k-1} \sum_{\un\in \Gk_{j,l}} \big(f_k(1, n_1, \dots, n_j)-f_k(z, n_1, \dots, n_j)\big).
\end{equation}
\end{lemma}

Having expressed $X_k(z)$ in terms of $X_0(z)$, for $k=1,2, \dots, g$, we can use \eqref{eqn:queuelengthevolution2} to find $X_c(z)$. Following the same steps as for deriving Lemma \ref{thm:xkx0}, we get the lemma below.

\begin{lemma}\label{lemma:xkz}
For $k=g+1,\dots,c$,
\begin{align*}
&X_k(z) = X_0(z)\frac{Y(\overbrace{z, \dots, z}^{k},1,\dots,1)}{z^g} \\
&+
\sum_{l=0}^{g-1}q_l \sum_{j=0}^{g-1} \sum_{\un\in \G^{(g)}_{j,l}} \big(h_g(1, \overbrace{z,\dots,z}^{k-g},\overbrace{1,\dots,1}^{c-k},n_1, \dots, n_j)-h_g(z, \overbrace{z,\dots,z}^{k-g},\overbrace{1,\dots,1}^{c-k}, n_1, \dots, n_j)\big)
\end{align*}
with $h_m(z, y_{m+1},\dots,y_c,n_1, \dots, n_j)$ defined as
\begin{equation}
\frac{\partial}{\partial y_1^{n_1}\partial y_2^{n_2}\dots \partial y_j^{n_j}} \,
 \frac{Y(y_1,y_2,\dots,y_j,\overbrace{z,\dots,z}^{m-j},y_{m+1},\dots,y_c)}{n_1!n_2!\dots n_j!\,z^{m-j}}\Big|_{y_1=\dots=y_j=0}.
\label{eqn:hg}
\end{equation}
\end{lemma}
If the system is in steady state, it should hold that $X_{c,n}\equaldist X_{0,n}$, so we can equate $X_{c}(z)$ and $X_{0}(z)$ and solve for $X_0(z)$ to arrive at the following result.
\begin{theorem}\label{thm:x0}
\begin{align}
X_0(z) = \frac{\sum_{l=0}^{g-1}q_l \zeta_l(z) }{z^g-Y(z, \dots, z)}, 
\label{maineqqq}
\end{align}
with
\begin{equation}
\zeta_l(z)=\sum_{j=0}^{g-1} \sum_{\un\in \G_{j,l}} z^g\big(h_g(1, z,\dots,z,n_1, \dots, n_j)-h_g(z, z,\dots,z,n_1, \dots, n_j)\big).
\end{equation}
\end{theorem}

\begin{rem}
To help interpret \eqref{maineqqq} write
\begin{align}
z^g h_g(1, z,\dots,z,n_1, \dots, n_j)&=\E\left[z^{g+Y_{g+1}+\dots+Y_c}\mathbf{1}_{\{Y_1=n_1, \dots, Y_j=n_j\}}\right],\label{zetal2}\\
z^g h_g(z, z,\dots,z,n_1, \dots, n_j)&=\E\left[z^{j+Y_{j+1}+\dots+Y_c}\mathbf{1}_{\{Y_1=n_1, \dots, Y_j=n_j\}}\right],\label{zetal1}
\end{align}
and
\begin{align}
X_0(z) &= \frac{Y(z, \dots, z)}{z^g}X_0(z)+\frac{\sum_{l=0}^{g-1}q_l \zeta_l(z) }{z^g} \nonumber\\
&= \frac{Y(z, \dots, z)}{z^g}X_0(z)+\sum_{l=0}^{g-1}q_l \E\left[\left(z^{g+Y_{g+1}+\dots+Y_c}-z^{j+Y_{j+1}+\dots+Y_c}\right)\mathbf{1}_{\{Y_1=n_1, \dots, Y_j=n_j\}}\right].
\label{eqx0alt} 
\end{align}

The first expression on the right hand side of \eqref{eqx0alt} represents the relation $X_{0,n+1}=X_{0,n} + \sum_{i=1}^c Y_{i,n}-g$ in terms of generating functions. The last term in \eqref{eqx0alt} contains all ``correction terms'' for those summands where the queue becomes empty at time slot $j$, because then $X_{0,n+1}=X_{0,n} + \sum_{i=1}^j Y_{i,n}+ \sum_{i=g+1}^c Y_{i,n} - j$, which can be written as
\[
X_{0,n+1}=X_{0,n} + \sum_{i=1}^c Y_{i,n}-g + \Big(g+\sum_{i=g+1}^g Y_{i,n}-j-\sum_{i=j+1}^c Y_{i,n}\Big),
\]
the last part being the correction term. The function $\zeta_l(z)$ is composed of all correction terms that correspond to $X_{0,n}=l$. Each correction term consists of the term \eqref{zetal1}, representing the part that should be removed, and its replacement term \eqref{zetal2}.
\end{rem}

\begin{rem}
If the $Y_{i}$ are independent (but not necessarily identically distributed) random variables, Equation \eqref{maineqqq} can be simplified. In this case, $Y(y_1, \dots, y_c) = \prod_{i=1}^c Y_i(z)$, which can be substituted in \eqref{eqn:hg} to find
\begin{align}
&h_g(z, z,\dots,z,n_1, \dots, n_j)=z^{j-g}\,\P(Y_1=n_1, \dots, Y_j=n_j)\prod_{i=j+1}^c Y_i(z),\\
&h_g(1, z,\dots,z,n_1, \dots, n_j)=\P(Y_1=n_1, \dots, Y_j=n_j)\prod_{i=g+1}^c Y_i(z).
\label{eqn:hgindep}
\end{align}
Let the random variable $G$ denote the effective green time, as defined earlier, and notice that
\begin{equation}
\P(G=j) = \sum_{l=0}^{g-1}q_l \sum_{\un\in \G_{j,l}} \P(Y_1=n_1, \dots, Y_j=n_j) .
\label{eqn:PGj}
\end{equation}

\noindent From \eqref{maineqqq} we then obtain
\begin{align}
X_0(z) &= \frac{z^g \sum_{j=0}^{g-1} \P(G=j) \left(\prod_{i=g+1}^c Y_i(z) - z^{j-g}\,\prod_{i=j+1}^c Y_i(z)\right)}{z^g-\prod_{i=1}^c Y_i(z)}\nonumber\\
&= \frac{(1-\P(G=g)) z^g  \prod_{i=g+1}^c Y_i(z) - \sum_{j=0}^{g-1} \P(G=j)  z^{j}\,\prod_{i=j+1}^c Y_i(z)}{z^g-\prod_{i=1}^c Y_i(z)}.
\label{eqn:X0indep}
\end{align}

\noindent Using that $\P(G=0)=\P(X_0=0)$, $\P(G=g)=1-\P(X_{g-1}=0)$, and
$\P(G=j)=\P(X_j=0)-\P(X_{j-1}=0)$, for $j=1,\dots,g-1$,
\eqref{eqn:X0indep} can be written as
\begin{equation*}
X_0(z) = \frac{\sum_{j=1}^g\left(1-Y_j(z)/z\right)\P(X_{j-1}=0) z^j\prod_{i=j+1}^c Y_i(z)}{z^g-\prod_{i=1}^c Y_i(z)},\label{eqn:X0independent}
\end{equation*}
where the unknown probabilities $\P(X_i=0)$ are still to be determined.
\end{rem}

\begin{rem}
When $Y_i$ are i.i.d. random variables, with pgf $Y(z)$, this expression simplifies further to the classical result \cite{darroch,fctlsolo}
\begin{equation*}
X_0(z) = \frac{\left(1-Y(z)/z\right)\sum_{j=1}^g \P(X_{j-1}=0)z^j Y(z)^{c-j}}{z^g-Y(z)^c}.\label{eqn:X0independent}
\end{equation*}
\end{rem}

Let us now proceed with the expression in \eqref{maineqqq} and show how this generating function can be converted into algorithms for the performance analysis of the generalized FCTL queue. First note that there are still $g$ unknowns $q_0,\ldots,q_{g-1}$ in \eqref{maineqqq}, which can be found using a classical reasoning. With Rouch\'{e}'s theorem \cite{adan2006application}, it can be shown that the denominator of (\ref{maineqqq}) has $g$ zeros on or within the unit circle $|z|\leq 1$. Since a pgf is analytic and
well-defined in $|z|\leq 1$, the numerator of $X_0(z)$ should vanish at each of the zeros. This gives
$g$ equations. One of the zeros equals 1, and leads to a trivial equation.
However, the normalization condition $X_0(1)=1$ provides an
additional equation. Using l'H\^{o}pital's rule, this condition is
found to be
\begin{equation}\label{normcondfixedfctlq}
\sum_{l=0}^{g-1}q_l\sum_{j=0}^{g-1} \sum_{\un\in \G_{j,l}} \Big(g-j-\sum_{i=j+1}^g\E[Y_i \,|\, A_j]\Big)\P(A_j)=g-\sum_{i=1}^c\E[Y_i],
\end{equation}
where we have defined $A_j$ as the event $\{Y_1=n_1,\dots,Y_j=n_j\}$ for compactness. From this definition it follows that $\E[Y_i \,|\, A_0]=\E[Y_i]$ and $\P(A_0)=1$.
We can write \eqref{normcondfixedfctlq} as
\[
\sum_{l=0}^{g-1}q_l b_l = \eta.
\]
Note that $b_l=\zeta_l'(1)$, with $\zeta_l(z)$ as defined in \eqref{maineqqq}.
Denote the $g$ roots of $z^g=Y(z, \dots, z)$ on and within the unit circle by $z_0=1,z_1,\ldots,z_{g-1}$.
The $g$ unknowns $q_0,\ldots,q_{g-1}$ then follow from solving the set of linear equations
\begin{equation}\label{fctlseteqnsss}
\begin{pmatrix}
  b_0 & b_1 & b_2 &\ldots & b_{g-1} \\
  \zeta_0(z_1)  & \zeta_1(z_1) & \zeta_2(z_1) &\ldots & \zeta_{g-1}(z_1) \\
    \zeta_0(z_2)  & \zeta_1(z_2) & \zeta_2(z_2) &\ldots & \zeta_{g-1}(z_2) \\
  \vdots & \vdots & \vdots &\vdots & \vdots \\
  \zeta_0(z_{g-1})  & \zeta_1(z_{g-1}) & \zeta_2(z_{g-1}) &\ldots & \zeta_{g-1}(z_{g-1})
\end{pmatrix}
\left(%
\begin{array}{c}
  q_0 \\
  q_1 \\
  q_2 \\
  \vdots \\
  q_{g-1} \\
\end{array}\right)=\left(\begin{array}{c}
  \eta \\
  0 \\
  0 \\
  \vdots \\
  0 \\
\end{array}%
\right).
\end{equation}
To make this work, we need an efficient numerical algorithm for finding the roots $z_1,\ldots,z_{g-1}$. Many root-finding algorithms exist, for instance based on successive substitution \cite{janssen2005analytic} or the Lagrange inversion theorem \cite{janssen2008back}. In the appendix we present a new scheme (Algorithm \ref{alg:roots}) based on the truncation of infinite series in combination with a root-finding procedure for polynomial functions.

We also need to deal with the sets $\Sk_{j,l}$ that may have an infinite number of elements. For this, notice that $\Gk_{j,l}$ has a finite cardinality. Finding all elements in $\Gk_{j,l}$ is rather straightforward, because they can be found by a simple enumeration. When $l=0$, the effective green period is zero, meaning that $j$ must be zero as well. As a consequence, all $\Gk_{j,0}$ are empty sets, except for $\Gk_{0,0}$, which is a set containing one element: the empty vector $(\,)$. For $0<l<k$, the minimum effective green time is at least $l$. Fixing $l$ between 1 and $k-1$, and $j$ between $l$ and $k-1$, we note that all elements in $\Gk_{j,l}$ must satisfy the following conditions:
\begin{align*}
l+n_1-1&>0, \\
l+n_1+n_2-2&>0,\\
&\ \,\vdots\\
l+n_1+\dots+n_{j-1}-(j-1)&>0,\\
l+n_1+\dots+n_{j}-j&=0.
\end{align*}
Since $n_1, \dots, n_j$ are nonnegative integers, these conditions lead to an easy enumeration of all elements $\un$ in $\Gk_{j,l}$. See Algorithm \ref{alg:Gjl} in the appendix for efficiently calculating all nonempty sets $\Gk_{j,l}$.

Taken everything together, we can evaluate the generating function $X_0(z)$ as described in Algorithm \ref{alg:x0pgf}.

\begin{algorithm}%
\caption{Computing $X_0(z)$}%
\label{alg:x0pgf}%
\begin{algorithmic}[1]%
\State Input $g$, $c$, and $Y(z_1, z_2, \dots, z_c)$
\State Set $\E[Y]=\frac{\text{d}}{\text{d}z}Y(z, z, \dots, z)\big|_{z=1}$
\State Generate all nonempty sets $\Gk_{j,l}$ using Algorithm \ref{alg:Gjl}
\State Determine all roots $z_1, z_2, \dots, z_{g-1}$ inside the unit circle using Algorithm \ref{alg:roots}
\For{$l=0$ to $g-1$}
  \State Set $\zeta_l(z) = 0$
  \For{$j=l$ to $g-1$}
    \ForAll{$\un \in \Gk_{j,l}$}
        \State Set $(n_1, n_2, \dots, n_j) = \un$
        \State Compute $h_g(z, y_{g+1},\dots,y_c,n_1, \dots, n_j)$ using Equation \eqref{eqn:hg}
        \State Set $\zeta_l(z) = \zeta_l(z) + z^g\big(h_g(1, z,\dots,z,n_1, \dots, n_j)-h_g(z, z,\dots,z,n_1, \dots, n_j)\big)$
    \EndFor
  \EndFor
  \State Set $b_l = \zeta_l'(1)$
\EndFor
\State Set $\eta=g-\E[Y]$
\State Find $q_1, q_2, \dots, q_{g-1}$ by solving system \eqref{fctlseteqnsss}
\State Compute $X_0(z)$ using Equation \eqref{maineqqq}
\State \Return $X_0(z)$
\end{algorithmic}%
\end{algorithm}%

\subsection{Queue length distribution}\label{s43}
Now that we can calculate $X_0(z)$, we proceed to use Algorithm \ref{alg:x0pgf} to obtain more information about the queue length. Notice that with Lemmas \ref{thm:xkx0} and~\ref{lemma:xkz} and Algorithm \ref{alg:x0pgf} we can evaluate all $X_k(z)$ for $k=1,\dots,c$. While these are all pgfs at the end of specific time slots, the pgf of the queue length at an arbitrary time slot $\bar{X}(z)$ follows from
\begin{equation}\label{XzAvg}
\bar{X}(z) = \frac1c\sum_{i=1}^c X_i(z).
\end{equation}
The queue-length distribution can then be calculated using a standard inverse theorem such as Algorithm \ref{alg:pgfinversion} in the appendix. We have used Algorithm \ref{alg:pgfinversion} with \eqref{XzAvg} to generate the results in Table \ref{tbl:queuelengthprobsExample1}. In the appendix we also present some compact expressions for the mean queue length at the end of each time slot.

\section{Network settings}\label{sec:networks}

With the algorithms developed in Section \ref{sec:main} to analyze the output of an intersection with correlated input, we now extend the scope in order to deal with a network of intersections. In Subsection \ref{dec} we describe the decomposition approach, a heuristic method to combine multiple isolated intersections into a network model, and in Subsection \ref{sect:networkexamples} we demonstrate the algorithm for two network scenarios.

\subsection{Decomposition approach}\label{dec}

Let us first quantify the output process of the generalized FCTL queue in Section \ref{sec:main}. Define $O_{i,n}$ as the output in time slot $i$ ($i=1,\dots,g$) in cycle $n$, so that
\begin{equation}
O_{i,n}=\begin{cases}
1 & \qquad \text{ if }X_{i-1,n}>0,\\
Y_{i,n} & \qquad \text{ if }X_{i-1,n}=0.
\end{cases}
\label{eqn:marginaloutput}
\end{equation}
Let $O_n(z_1,\dots,z_g)$ denote the pgf of the joint output $(O_{1,n},\dots,O_{g,n})$ and $G_n$ denote the effective green time in cycle $n$. Note that the joint output vector $(O_{1,n},\dots,O_{g,n})$, given that $G_n=j$, equals
\[
(\overbrace{1, 1, \dots, 1}^{j},Y_{j+1,n}, \dots, Y_{g,n}),\qquad j=1,2,\dots,g.
\]
We will again omit the subscript $n$ denoting the cycle number. If $(n_1,\dots,n_{j})$ is an arbitrary vector in $\G_{j,l}$, it holds that $l+n_1+\dots+n_j=j$. The pgf of the joint output during this cycle follows from summing over all possible $l, j$ and $(n_1,\dots,n_{j}) \in \G_{j,l}$:
\begin{align}
O&(z_1,\dots,z_g)=
\sum_{l=0}^\infty \sum_{n_1=0}^\infty\dots\sum_{n_g=0}^\infty \big(1-T(l,n_1,\dots,n_{g-1})\big) \Big(\prod_{k=1}^g z_k\Big)\,
\P(X_{0}=l,Y_{1}=n_1,\dots,Y_{g}=n_g)\nonumber\\
&+
\sum_{l=0}^{g-1}\sum_{j=0}^{g-1}
\sum_{(n_1,\dots,n_{j})\in \G_{j,l}}\sum_{n_{j+1}=0}^\infty\dots\sum_{n_g=0}^\infty\Big(\prod_{k=1}^j z_k\Big) \Big(\prod_{k=j+1}^g z_k^{n_k}\Big)\P(X_{0}=l,Y_{1}=n_1,\dots,Y_{g}=n_g).
\end{align}
This then gives
\begin{align}
O&(z_1,\dots,z_g)=\nonumber\\&=
\Big(1-
\sum_{l=0}^{g-1}\sum_{j=0}^{g-1}
\sum_{(n_1,\dots,n_{j})\in \G_{j,l}}\P(X_{0}=l,Y_{1}=n_1,\dots,Y_{j}=n_j)\Big)\prod_{k=1}^g z_k\nonumber\\
&+
\sum_{l=0}^{g-1}\sum_{j=0}^{g-1}
\sum_{(n_1,\dots,n_{j})\in \G_{j,l}}\sum_{n_{j+1}=0}^\infty\dots\sum_{n_g=0}^\infty\Big(\prod_{k=1}^j z_k\Big)\Big(\prod_{k=j+1}^g z_k^{n_k}\Big)\P(X_{0}=l,Y_{1}=n_1,\dots,Y_{g}=n_g)\nonumber\\
&=
\Big(1-
\sum_{l=0}^{g-1}\sum_{j=0}^{g-1}
\sum_{\un\in \G_{j,l}}q_l 
\P(A_j)\Big)\,\prod_{k=1}^g z_k\nonumber \\
&+
\sum_{l=0}^{g-1} q_l \,\sum_{j=0}^{g-1}
\sum_{\un\in \G_{j,l}} \Big(\prod_{k=1}^j z_k\Big) h_{j}(1, z_{j+1},\dots,z_{g},\overbrace{1,\dots,1}^{c-g},n_1,\dots,n_j)
\end{align}
with $h_j$ as defined in \eqref{eqn:hg}. Note that $h_0(1, z_{1},\dots,z_{g},1,\dots,1)=Y(z_1,\dots,z_g,1,\dots,1)$ and $\P(A_j)=h_j(1, 1,\dots,1, n_1,\dots,n_j)$.
The complete program to compute $O(z_1, z_2, \dots, z_g)$ is then given in Algorithm \ref{alg:outputpgf}.

\begin{algorithm}[!ht]%
\caption{Computing $O(z_1, z_2, \dots, z_g)$}%
\label{alg:outputpgf}%
\begin{algorithmic}[1]%
\State Input $g$, $c$, and $Y(z_1, z_2, \dots, z_c)$
\State Find $q_1, q_2, \dots, q_{g-1}$ using Algorithm \ref{alg:x0pgf}
\State Set $S=0, T=0$
\For{$l=0$ to $g-1$}
  \For{$j=l$ to $g-1$}
    \ForAll{$\un \in \Gk_{j,l}$}
        \State Set $(n_1, n_2, \dots, n_j) = \un$
        \State Compute $h_j(z, y_{j+1},\dots,y_c,n_1, \dots, n_j)$ using Equation \eqref{eqn:hg}
        \State Set $S=S+q_l h_j(1, 1,\dots,1,n_1, \dots, n_j)$
        \State Set $T=T+q_l\left(\prod_{k=1}^j z_k\right) h_{j}(1, z_{j+1},\dots,z_{g},1,\dots,1,n_1,\dots,n_j),$
    \EndFor
  \EndFor
\EndFor
\State Set $O(z_1, z_2, \dots, z_g)=(1-S)\prod_{k=1}^g z_k+T$
\State \Return $O(z_1, z_2, \dots, z_g)$
\end{algorithmic}%
\end{algorithm}%

\begin{rem}
It is evident that correlation in input carries over to correlation in output. Less obvious is that even with independent input, the output processes in consecutive time slots are correlated. In the case of independent $Y_1,\dots,Y_c$, the joint output pgf reduces to
\begin{align}
O(z_1, z_2, \dots, z_g)
=&\P(G=0) \prod_{i=1}^g Y_i(z_i)+ \sum_{j=1}^{g-1}\P(G=j)\prod_{i=1}^{j} z_i\, \prod_{i=j+1}^g Y_i(z_i)\nonumber\\
 & + \P(G=g) \prod_{i=1}^g z_i.
\label{eqn:jointoutputindependent}
\end{align}
The correlation follows from the fact that this is not the product of the marginal pgfs $O_i(z_i)$ for $i=1,\dots, g$.
\end{rem}

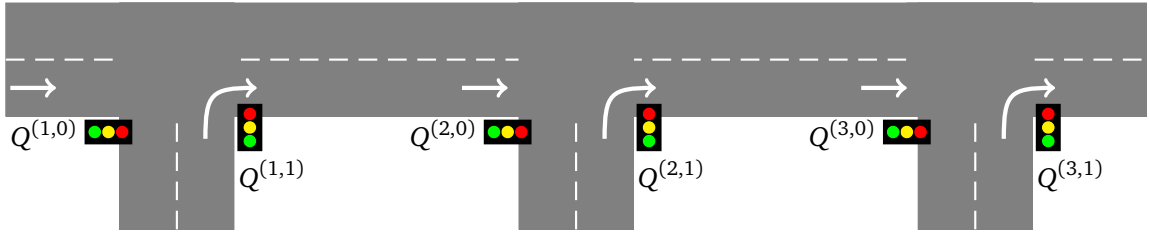
\begin{figure}[ht]
\begin{center}
\begin{tikzpicture}[scale=0.3]
\draw[gray,fill=gray](5,0) rectangle (10,10);
\draw[thick,white,dash pattern=on 7 off 4](7.5,0) to (7.5,5);
\draw[gray,fill=gray](22.5,0) rectangle (27.5,10);
\draw[thick,white,dash pattern=on 7 off 4](25,0) to (25,5);
\draw[gray,fill=gray](40,0) rectangle (45,10);
\draw[thick,white,dash pattern=on 7 off 4](42.5,0) to (42.5,5);
\draw[gray,fill=gray](0,5) rectangle (50,10);
\draw[thick,white,dash pattern=on 7 off 4](0,7.5) to (50,7.5);
\draw[gray,fill=gray](5,5) rectangle (10,10);
\draw[gray,fill=gray](22.5,5) rectangle (27.5,10);
\draw[gray,fill=gray](39.5,5) rectangle (45,10);
\draw[ultra thick,white,->](0.2,6.25) to (2.2,6.25);
\draw[ultra thick,white,->](20,6.25) to (22,6.25);
\draw[ultra thick,white,->](37.5,6.25) to (39.5,6.25);
\draw[ultra thick,white,->](8.75,4.0) to [out=90, in=180, distance=2cm] (11.0,6.25);
\draw[ultra thick,white,->](26.25,4.0) to [out=90, in=180, distance=2cm] (28.5,6.25);
\draw[ultra thick,white,->](43.75,4.0) to [out=90, in=180, distance=2cm] (46.0,6.25);
\trafficsignal[-90]{(3.5,4.8)}
\trafficsignal[-90]{(21,4.8)}
\trafficsignal[-90]{(38.5,4.8)}
\node[below] at (1.7,5.3) {$Q^{(1,0)}$};
\node[below] at (19.2,5.3) {$Q^{(2,0)}$};
\node[below] at (36.7,5.3) {$Q^{(3,0)}$};
\trafficsignal{(10.2,3.5)}
\trafficsignal{(27.7,3.5)}
\trafficsignal{(45.2,3.5)}
\node[below] at (11.7,3.5) {$Q^{(1,1)}$};
\node[below] at (29.2,3.5) {$Q^{(2,1)}$};
\node[below] at (46.7,3.5) {$Q^{(3,1)}$};
\end{tikzpicture}
\end{center}
\caption{The network considered in this section.}
\label{fig:networksidetraffic}
\end{figure}

To analyze a network of fixed-cycle intersections, we decompose the network into isolated generalized FCTL queues. We exclude situations where the output of a queue could become the input of that same queue through some cyclic path in the network, meaning that the network can be represented as a directed acyclic graph where nodes represent queues. Nodes without parents can be modeled as standard FCTL queues with independent external arrivals. Figure \ref{fig:networksidetraffic} shows an example of three intersections with one main traffic flow, with queues labeled $Q^{(1,0)}, Q^{(2,0)}$, and $Q^{(3,0)}$, and three minor flows $Q^{(1,1)}, Q^{(2,1)}$, and $Q^{(3,1)}$ representing side traffic.

This network can be decomposed into six queues. Queues $Q^{(1,0)}, Q^{(1,1)}, Q^{(2,1)}$, and $Q^{(3,1)}$ are standard FCTL queues with Poisson input. The analysis in \cite{fctlsolo} can be used to obtain queue length distributions and Equation \eqref{eqn:jointoutputindependent} gives the joint output of each of these queues. Following the main flow downstream, we first determine the input of $Q^{(2,0)}$, which is the superposition of the delayed output processes of $Q^{(1,0)}$ and $Q^{(1,1)}$. The analysis in Section \ref{sec:main} gives the queue length pgf, and the results from this section give the joint output. Algorithm \ref{alg:network} provides a detailed description of this method, which can be used to find approximations for all queue-length pgfs in the network.

\begin{algorithm}[!ht]%
\caption{Network analysis}%
\label{alg:network}%
\begin{algorithmic}[1]%
\State Input the number of intersections $n$
\State Input $g^{(1,0)}$, $c$, and $Y^{(1,0)}(z_1, z_2, \dots, z_c)$  for $Q^{(1,0)}$
\State Compute $X^{(1,0)}_0(z)$ using Algorithm \ref{alg:x0pgf}%
\State Compute $O^{(1,0)}(z_1, z_2, \dots, z_{g^{(1,0)}})$ using Algorithm \ref{alg:outputpgf}%
\For{$i=2$ to $n$}
   \State Input $g^{(i-1,1)}$ and $Y^{(i-1,1)}(z_1, z_2, \dots, z_c)$  for $Q^{(i-1,1)}$
   \State Input $\delta$, the number of red time slots preceding the first green time slot at $Q^{(i-1,1)}$
   \State Compute $O^{(i-1,1)}(z_1, z_2, \dots, z_{g^{(i-1,1)}})$ using Algorithm \ref{alg:outputpgf}%
    \State Set $Y^{(i,0)}(z_1, z_2, \dots, z_c) = O^{(i-1,0)}(z_{\text{mod}(d, c)+1}, z_{\text{mod}(1+d, c)+1}, \dots, z_{\text{mod}({g^{(i-1,0)}}-1+d, c)+1})$
    \State $\hspace*{0.3cm}\times O^{(i-1,1)}(z_{\text{mod}(\delta+d, c)+1}, z_{\text{mod}(1+\delta+d, c)+1}, \dots, z_{\text{mod}({g^{(i-1,1)}}-1+\delta+d, c)+1})$
\State Compute $X^{(i,0)}_0(z)$ using Algorithm \ref{alg:x0pgf}%
\State Compute $O^{(i,0)}(z_1, z_2, \dots, z_{g^{(i,0)}})$ using Algorithm \ref{alg:outputpgf}%
\EndFor
\end{algorithmic}%
\end{algorithm}%

\subsection{Two example networks}\label{sect:networkexamples}

Consider an idealized network of ten connected generalized FCTL queues, $Q^{(1)}, \dots, Q^{(10)}$, each with a fixed cycle of $g=10$ green and $r=10$ red time slots and a fixed travel time between two consecutive intersections of $d$ time slots. The setting is similar to that in Figure~\ref{fig:networksidetraffic} without the side traffic. Since all intersections have the same fixed-cycle plan of 10 green and 10 red time slots, the case $d=0$ corresponds to a perfect green wave where vehicles arrive exactly at the moment that the signal turns green, and no queues will build up. In the case $d=5$, the vehicles departing during the first five time slots of $Q^{(i)}$, will arrive at $Q^{(i+1)}$ during time slots 6--10. These vehicles will pass without delay, except when there is a queue of more than five vehicles waiting at the beginning of the cycle. The vehicles departing in time slots 6--10 from $Q^{(i)}$, will arrive at $Q^{(i+1)}$ during time slots 11--15, which means that they arrive during the red period and will be delayed until (at least) the next green period.

Only $Q^{(1)}$ has external arrivals, assumed to arrive according to a Poisson process with rate $\lambda=0.45$ per time slot. The \emph{occupation rate} $\rho:=\lambda c/g=0.9$ is close to one, indicating that the intersections are operating close to their maximum capacity. Appendix \ref{appendix:input} provides more details about the input settings of this example, and how to determine the arrival processes at the other intersections.

\begin{figure}[ht]
\begin{center}

\begin{tabular}{cc}
\hspace*{-6.5cm}$\E[\bar{X}^{(i)}]$ & \hspace*{-6.5cm}$\E[\bar{X}^{(i)}]$ \\
\includegraphics[width=0.45\textwidth]{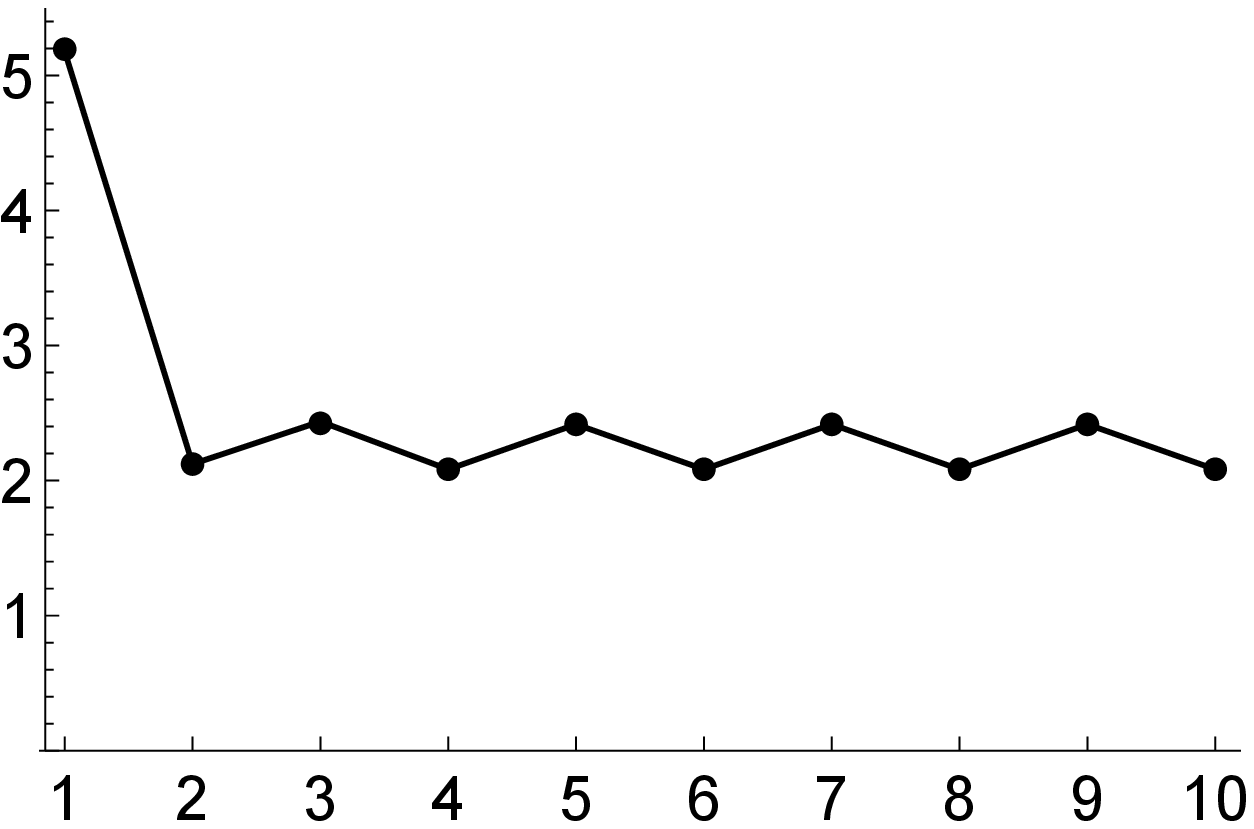} \raisebox{3mm}{$i$}& \includegraphics[width=0.45\textwidth]{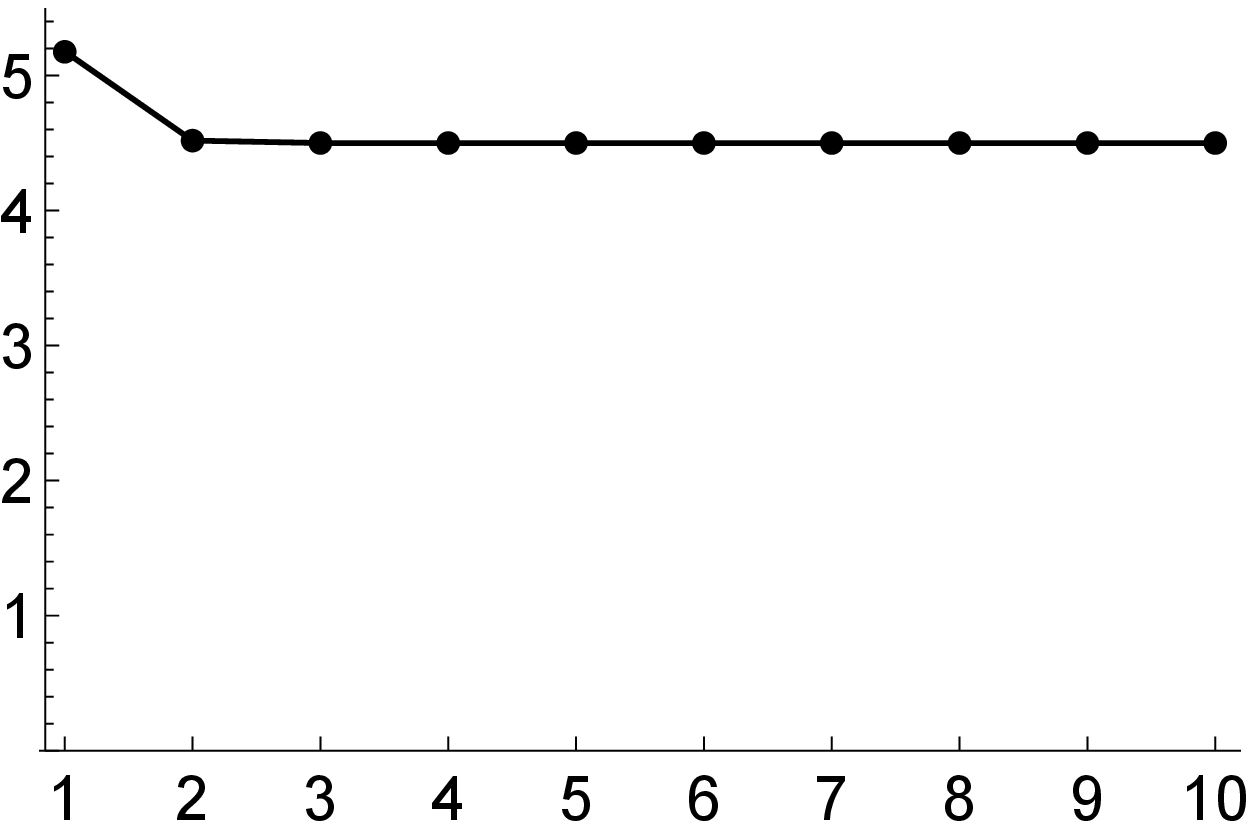} \raisebox{3mm}{$i$}\\
(a) 
$d=5$&(b) 
$d=10$\\[1ex]
\end{tabular}
\end{center}
\caption{Overall mean queue lengths in the network without side traffic in Section \ref{sect:networkexamples}. The dots represent the network simulation, the solid lines are obtained via the decomposition method. Figures \ref{fig:networknosidetraffic2}(a) and \ref{fig:networknosidetraffic2}(b) depict the expected overall mean queue length at intersection $i$ ($i=1,2,\dots,10$), $\E[\bar{X}^{(i)}] := \sum_{k=1}^c \E[X^{(i)}_k]/c$  for the cases $d=5$ and $d=10$ respectively. }
\label{fig:networknosidetraffic2}
\end{figure}

Observe that in Figures \ref{fig:networknosidetraffic2}(a) and \ref{fig:networknosidetraffic2}(b) the simulated values (the black dots) are indistinguishable from those computed using our generalized FCTL analysis, based on the decomposition approach (the black solid lines). Clearly, $d=5$ results in a better traffic flow, with smaller mean queue lengths, than $d=10$.

We now extend the previous example by merging traffic in the main flow with side traffic from minor roads, travelling towards the same destination as illustrated in Figure \ref{fig:networksidetraffic}. Each intersection has two flows, the main flow and the minor flow with side traffic. 
Note that the arrival patterns of queues in the main flows are similar to those in Section~\ref{s32}, with one large batch of vehicles arriving from the upstream main flow, and one smaller batch from the upstream minor flow. Each of the batches may be followed by vehicles arriving in the free flow. The specific input settings in this example can be found in Appendix \ref{appendix:input}.

\begin{figure}[!ht]
\begin{center}
\begin{tabular}{cc}
\hspace*{-6.5cm}$\E[{X}_i^{(2,0)}]$ & \hspace*{-6.5cm}$\E[{X}_i^{(8,0)}]$ \\[1mm]
\includegraphics[width=0.45\textwidth]{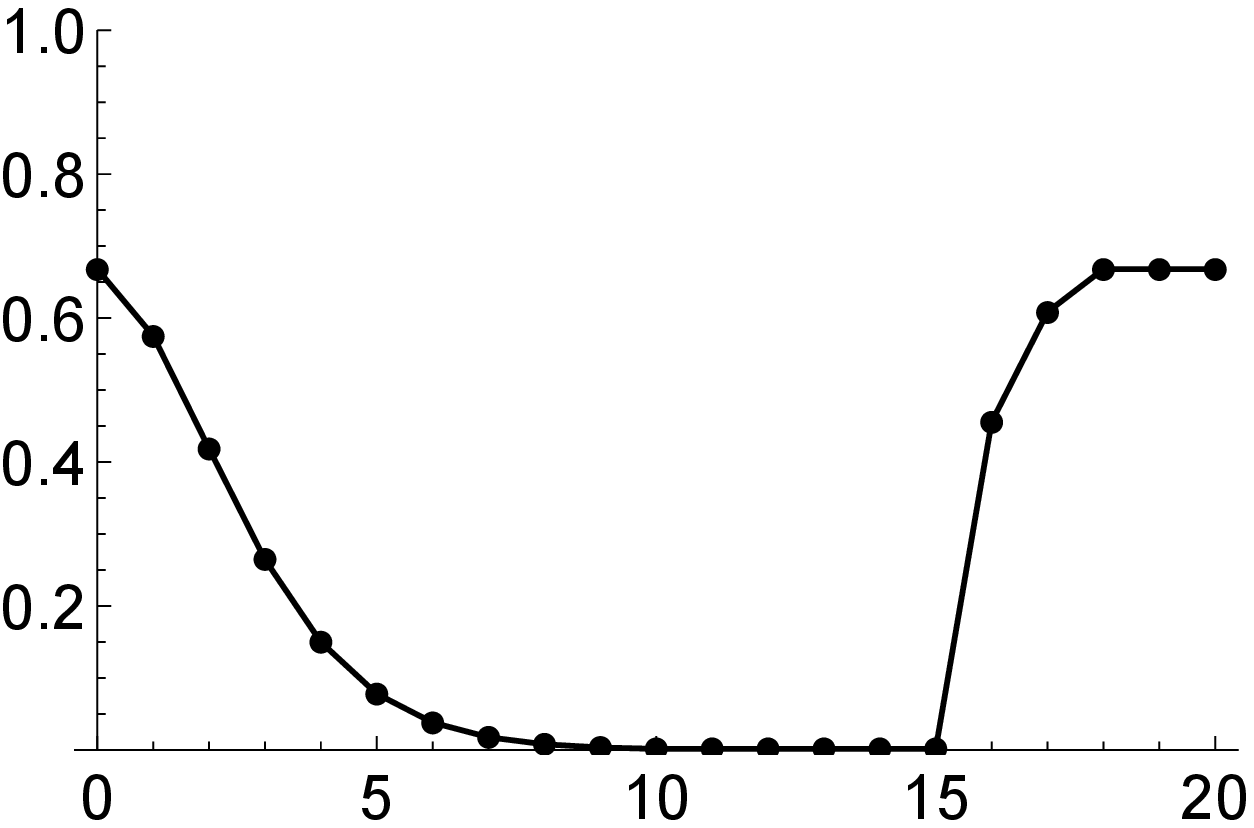} \raisebox{3mm}{$i$}& \includegraphics[width=0.45\textwidth]{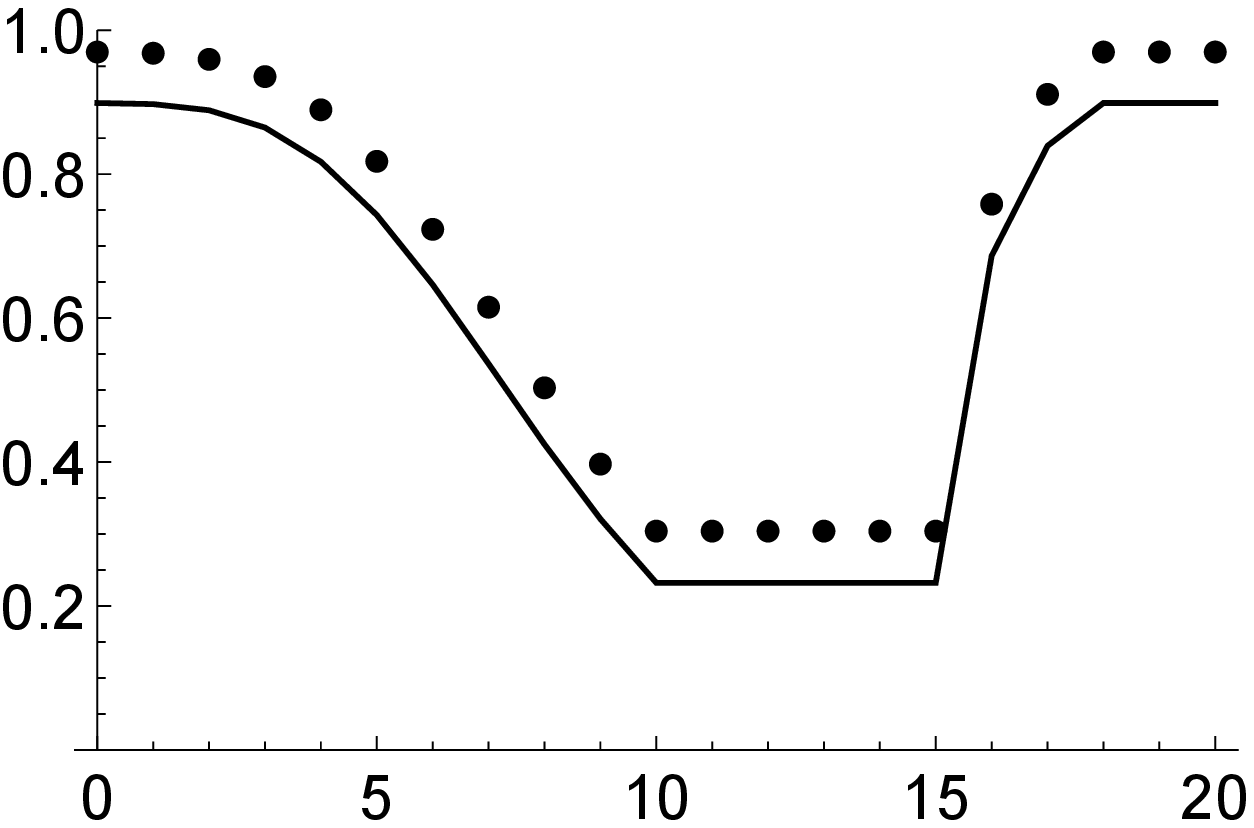} \raisebox{3mm}{$i$}\\
(a) $Q^{(2,0)}$ with $d=0$& (b) $Q^{(8,0)}$ with $d=0$\\[1ex]
\hspace*{-6.5cm}$\E[{X}_i^{(2,0)}]$ & \hspace*{-6.5cm}$\E[{X}_i^{(8,0)}]$ \\[1mm]
\includegraphics[width=0.45\textwidth]{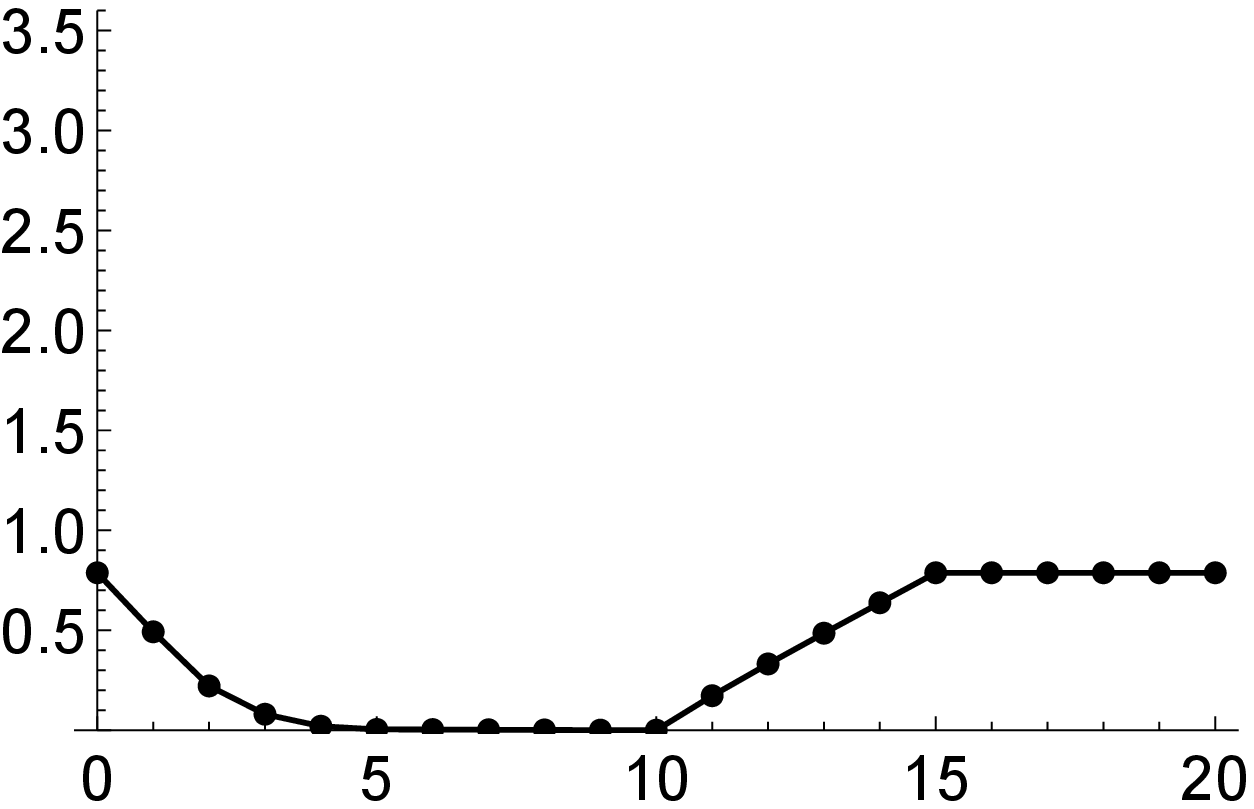} \raisebox{3mm}{$i$}& \includegraphics[width=0.45\textwidth]{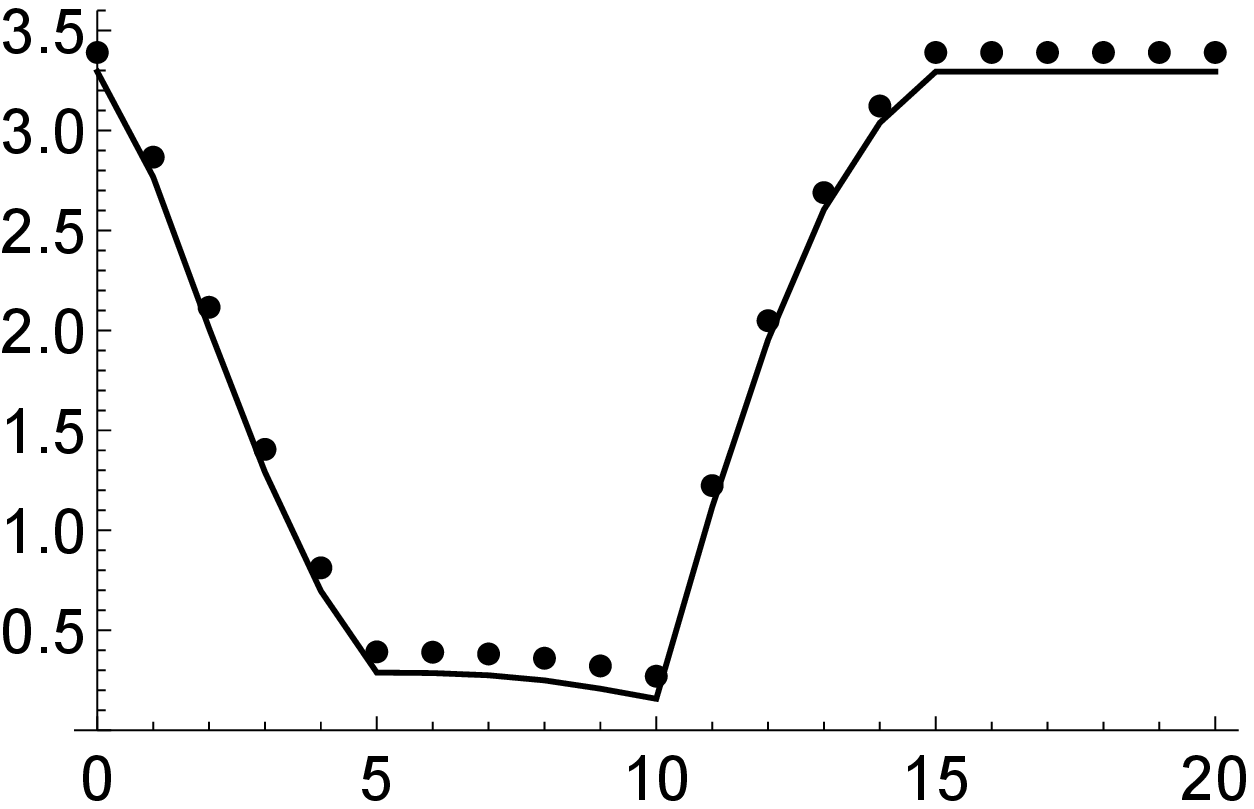} \raisebox{3mm}{$i$} \\
(c) $Q^{(2,0)}$ with $d=5$& (d) $Q^{(8,0)}$ with $d=5$\\[1ex]
\end{tabular}
\end{center}
\caption{Mean queue lengths for intersections 2 and 8 in the network with side traffic in Section \ref{sect:networkexamples}. The dots represent the network simulation, the solid lines are obtained via the decomposition method.}
\label{fig:networksidetraffic2}
\end{figure}

We have implemented this network model using Algorithm \ref{alg:network}, with the purpose of comparing the mean queue lengths of the main flows at each of the 10 intersections, but now for the cases $d=0$ and $d=5$. Due to the inflow of side traffic, the case $d=0$ will no longer result in empty queues. The large batch of (maximally) ten vehicles will experience no delay \emph{only if} no inflow from side traffic took place during the previous cycle. With $d=5$, the small batch of (maximally) three vehicles will arrive exactly at the moment that the traffic signal of the main flow turns green. As a consequence, these settings are favorable for the small batch, but also the first five vehicles from the large batch might benefit. The results in Figures \ref{fig:networksidetraffic2}(a)--(d) and in Table~\ref{tbl:networksidetraffic3} indicate that for intersections 1--7, the generalized FCTL queue gives extremely accurate approximations for the original network model, while for intersections 8--10, when the occupation rate $\rho$ exceeds 0.7, the algorithmic method is less accurate but still reasonable.

\begin{table}[!ht]
\begin{center}
\begin{tabular}{l|cccccccccc}
\hline
\multicolumn{11}{c}{Mean queue length $\bar{X}^{(i,0)}$ for $d=0$}\\
\hline
Intersection & 1 & 2 & 3 & 4 & 5 & 6 & 7 & 8 & 9 & 10\\
Load $\rho$ &               0.3 & 0.37 & 0.43 & 0.5 & 0.57 & 0.63 & 0.7 & 0.77 & 0.83 & 0.9 \\
\hline
Simulation &           0.493 & 0.232 & 0.260 & 0.293 & 0.335 & 0.395 & 0.491 & 0.665 & 1.046 & 2.192 \\
Decomposition  &     0.493 & 0.231 & 0.260 & 0.292 & 0.333 & 0.386 & 0.464 & 0.588 & 0.810 & 1.323 \\
\hline
\end{tabular}

\bigskip
\begin{tabular}{l|cccccccccc}
\hline
\multicolumn{11}{c}{Mean queue length $\bar{X}^{(i,0)}$ for $d=5$}\\
\hline
Intersection & 1 & 2 & 3 & 4 & 5 & 6 & 7 & 8 & 9 & 10\\
Load $\rho$ &               0.3 & 0.37 & 0.43 & 0.5 & 0.57 & 0.63 & 0.7 & 0.77 & 0.83 & 0.9 \\
\hline
Simulation &            0.493 & 0.359 & 1.161 & 0.821 & 1.549 & 1.287 & 1.984 & 1.939 & 2.814 & 3.829 \\
Decomposition  &      0.493 & 0.359 & 1.159 & 0.819 & 1.534 & 1.273 & 1.920 & 1.835 & 2.478 & 2.858 \\
\hline
\end{tabular}
\end{center}
\caption{The mean queue lengths at each of the ten intersections, for the network with side traffic of Section \ref{sect:networkexamples}.}
\label{tbl:networksidetraffic3}
\end{table}

\section{Conclusions}\label{sec:conclusions}

Classical models for fixed-cycle intersections do not cover
scenarios where arriving traffic is correlated, as typically encountered in
traffic networks. Motivated by this observation, we introduced a
generalized fixed-cycle traffic model that can deal with an arbitrary arrival patterns such as platoons and non-homogeneous traffic intensities. Generating-function methods for the stationary queue-length distribution were used to develop an efficient numerical scheme for the performance analysis of an isolated intersection. We also proposed a network algorithm, which decomposes a network of intersections into separate generalized fixed-cycle network models, whose distributional output and input are matched according to the network layout. Extensive simulation experiments for a line network of multiple intersection support the network algorithm, and reveal only a slight loss of accuracy in heavy-traffic scenarios.  We expect
our network algorithm to work for many realistic network configurations, the study of which is an
interesting topic for further research.

\section*{Acknowledgments}

This work is supported by the NWO Gravitation Networks grant 024.002.003.  
The work of JvL is further supported by an NWO TOP-GO grant and by an ERC Starting Grant.

\newpage

\appendix

\section{Additional algorithms}
\begin{algorithm}[!ht]%
\caption{Finding the roots inside the unit circle of $z^g-Y(z,z,\dots,z)$}%
\label{alg:roots}%
\begin{algorithmic}[1]%
\State Input $g$ and $Y(z, z, \dots, z)$
\State Set $D(z) = z^g-Y(z,z,\dots,z)$
\State Set the truncation parameter $n=g + 90$
\State Set $R=\{\,\}$
\While{$|R| \neq g-1$}
    \State $n = n + 10$
    \State Numerically compute $\tilde{D}(z):=\sum_{k=0}^n a_k z^k$,
    \State \qquad the $n$-th order Taylor expansion of $D(z)$ about $z=0$
    \State Find the roots of $\tilde{D}(z)$ using a standard numerical polynomial root finding algorithm
    \State Denote these roots by $\tilde{z}_1, \tilde{z}_2, \dots, \tilde{z}_n$
    \State $R=\{\,\}$
    \For{$i=1$ to $n$}
        \State Numerically find a root $z_i$ for $D(z)$ using starting point $z=\tilde{z}_i$
        \If{$z_i \notin R$ \emph{and} $|z_i| < 1$}
        \State add $z_i$ to $R$
        \EndIf
    \EndFor
\EndWhile
\State \Return $R$
\end{algorithmic}%
\end{algorithm}%

\begin{algorithm}[!ht]%
\caption{Creating all nonempty sets $\Gk_{j,l}$}%
\label{alg:Gjl}%
\begin{algorithmic}[1]%
\State $G_{0,0}=\{(\,)\}$
\For{$l=1$ to $g-1$}
    \For{$j=l$ to $g-1$}
\State $G_{j,l} = \{\,\}$
\For{$n_1=\max(0, 2-l)$ to $j-l$}
\For{$n_2=\max(0, 3-l-n_1)$ to $j-l-n_1$}
\State $\vdots$
\For{$n_{j-2}=\max(0, j-1-l-n_1-\dots-n_{j-3})$ to $j-l-n_1-\dots-n_{j-3}$}
\State $n_{j-1}=\max(0, j-l-n_1-\dots-n_{j-2})$
\State $n_j=0$
\State $G_{j,l} = G_{j,l} \cup \{ (n_1,n_2,\dots,n_j) \}$
\EndFor
\State $\vdots$
\EndFor
\EndFor
\EndFor
\EndFor
\end{algorithmic}%
\end{algorithm}%

\begin{algorithm}[!ht]%
\caption{Determine $P(X=k)$ by inverting $X(z)$ \cite{abatewhittinversion,choudhurywhittinversion}}%
\label{alg:pgfinversion}%
\begin{algorithmic}[1]%
\State Input $X(z)$ and $k$
\State $l=1$
\State $\gamma=10$
\If{$k=0$}
  \State $p_k = X(10^{-\gamma})$
\Else
  \State $r=10^{-\gamma/(2lk)}$
  \State $p_k = \frac{1}{2lkr^{k}} \left(X(r)+(-1)^k X(-r)+2 \sum _{j=1}^{lk-1} \text{Re}\Big[X\left(r \exp\big({\frac{\pi  i j}{l k}}\big)\exp\big({-\frac{\pi  i j}{l}}\big) \right)\Big]\right)$
\EndIf
\State \Return $p_k$
\end{algorithmic}%
\end{algorithm}%

\section{Mean queue length}
Denote by $N(z)$ and $D(z)$ the numerator and denominator of $X_0(z)$ as given in \eqref{maineqqq}. By differentiation and L'H\^opital's rule, it can be shown that
\[
\E[X_0]=\frac{N''(1)-D''(1)}{2D'(1)},
\]
where we have used that $X_0(1)=1$ and hence $N'(1)=D'(1)$. Define
\[
Y^{(j)}:=\sum_{i=j+1}^g Y_i, \qquad Y^{(r)}=\sum_{i=g+1}^c Y_i, \qquad Y=\sum_{i=1}^c Y_i.
\]
It is readily checked that
\[
D'(1)=g-\E[Y], \qquad D''(1)=g(g-1)-\E[Y(Y-1)].
\]
For the numerator, we find
\begin{align*}
N'(1)&=\sum_{l=0}^{g-1}q_l\sum_{j=0}^{g-1} \sum_{\un\in \G_{j,l}}\left(
\E\left[(g+Y^{(r)})\,\mathbf{1}_{\{A_j\}}\right]
-
\E\left[(j+Y^{(j)}+Y^{(r)})\mathbf{1}_{\{A_j\}}\right]
\right) \\
&=\sum_{l=0}^{g-1}q_l\sum_{j=0}^{g-1} \sum_{\un\in \G_{j,l}} \left(g-j-\E[Y^{(j)} \,|\, A_j]\right)\P(A_j),
\end{align*}
where as before $A_j=\{Y_1=n_1,\dots,Y_j=n_j\}$. Taking the second derivative yields
\begin{align*}
N''(1)&=\sum_{l=0}^{g-1}q_l\sum_{j=0}^{g-1} \sum_{\un\in \G_{j,l}}\left(
\E\left[(g+Y^{(r)})^2\,\mathbf{1}_{\{A_j\}}\right]
-
\E\left[(j+Y^{(j)}+Y^{(r)})^2\mathbf{1}_{\{A_j\}}\right]
\right) -N'(1).
\end{align*}
Using $D''(1)=g^2-\E[Y^2]-D'(1)$ and $N'(1)=D'(1)$ gives the following result:

\begin{lemma}
\begin{align}
\E[X_0]&=\frac{
\E[Y^2]-g^2+\sum q_l \P(A_j)(g^2-j^2)}{2(g-\E[Y])}\nonumber\\
&-\sum_{l=0}^{g-1}\sum_{j=0}^{g-1} \sum_{\un\in \G_{j,l}} \frac{q_l}{2(g-\E[Y])}\P(A_j)\E\left[2(g-j)Y^{(r)}-2(j-Y^{(r)})Y^{(j)}-\big(Y^{(j)}\big)^2\,\big|\,A_j\right]. \label{expmean}
\end{align}
\end{lemma}

Observe that \eqref{expmean} depends on the $q_l$'s, $\E[Y_i]$ for $i=1,\dots,c$, and $\E[Y_iY_k \,|\,A_j]$ for $i=1,\dots,g$ and $k=1,\dots,c$.
The mean queue lengths $\E[X_k]$ for $k=1,2,\dots,c$ can be found by expressing $\E[X_k]$ in terms of $\E[X_{k-1}]$ or $\E[X_{k+1}]$. For queue lengths during the red period, we see that
\[
\E[X_k]=\E[X_{k+1}]-\E[Y_{k+1}], \qquad k=c-1,c-2,\dots,g,
\]
which means that these mean queue lengths can be expressed in $\E[X_c]$ (and hence $\E[X_0]$) by successive substitution.

During the green period, the expressions become slightly more involved, but it is possible to express $\E[X_k]$ in terms of $\E[X_{k-1}]$ as
\[
\E[X_k]=\E[X_{k-1}]+\E[Y_{k}]-1 - \sum_{l=0}^{k-1}q_l \sum_{j=0}^{k-1} \sum_{\un\in \G_{j,l}} \P(A_j)\left(\E[Y_{k+1}\,|\,A_j]-1\right), \ \text{ for }k=1,\dots,g-1.
\]

\section{Details of numerical examples}\label{appendix:input}

For completeness and reproducibility, we give a detailed overview of the input parameters of all numerical examples in this paper.

\paragraph{Example Subsection \ref{s32}.}

\begin{figure}[ht]
\begin{center}
\begin{tikzpicture}[scale=0.3]
\draw[gray,fill=gray](0,0) rectangle (50,2.5);
\draw[ultra thick,white,->](0.2,1.25) to (2,1.25);
\draw[ultra thick,white,->](46,1.25) to (47.8,1.25);
\draw[white,dash pattern=on 2 off 2,rounded corners](3,0.3) rectangle (6,2.2);
\draw[white,dash pattern=on 2 off 2,rounded corners](8.5,0.3) rectangle (16.5,2.2);
\draw[white,dash pattern=on 2 off 2,rounded corners](19.5,0.3) rectangle (32.5,2.2);
\draw[white,dash pattern=on 2 off 2,rounded corners](34.5,0.3) rectangle (45,2.2);
\trafficsignal{(48.4,-0.3)}
\car{(3.5,0.65)}
\car{(9.0,0.65)}
\car{(11.5,0.65)}
\car{(14.0,0.65)}
\car{(20.0,0.65)}
\car{(24.0,0.65)}
\car{(30.0,0.65)}
\car{(35.0,0.65)}
\car{(37.5,0.65)}
\car{(40.0,0.65)}
\car{(42.5,0.65)}
\node[below] at (3.2,0) {free flow};
\node[below] at (12.5,0) {platoon 2};
\node[below] at (26,0) {free flow};
\node[below] at (40,0) {platoon 1};
\end{tikzpicture}

(a)

\bigskip
\begin{tabular}{ccc}
\includegraphics[height=4cm]{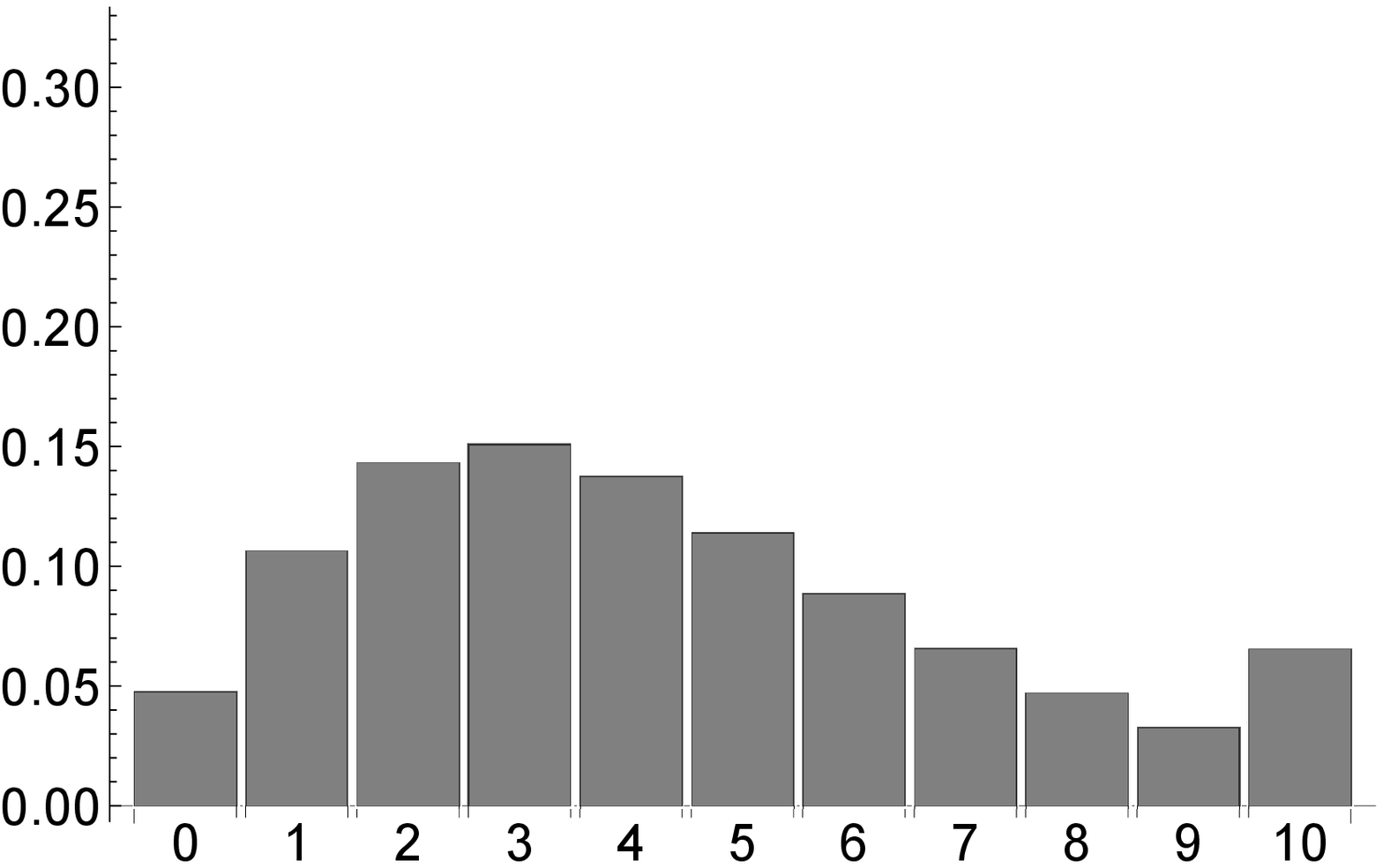}
\begin{picture}(0,0)
\put(-0.2,0.1){ $k$}
\put(-7,4.2){ $\P(B_1=k)$}
\end{picture}
&\hspace*{2cm}&\includegraphics[height=4cm]{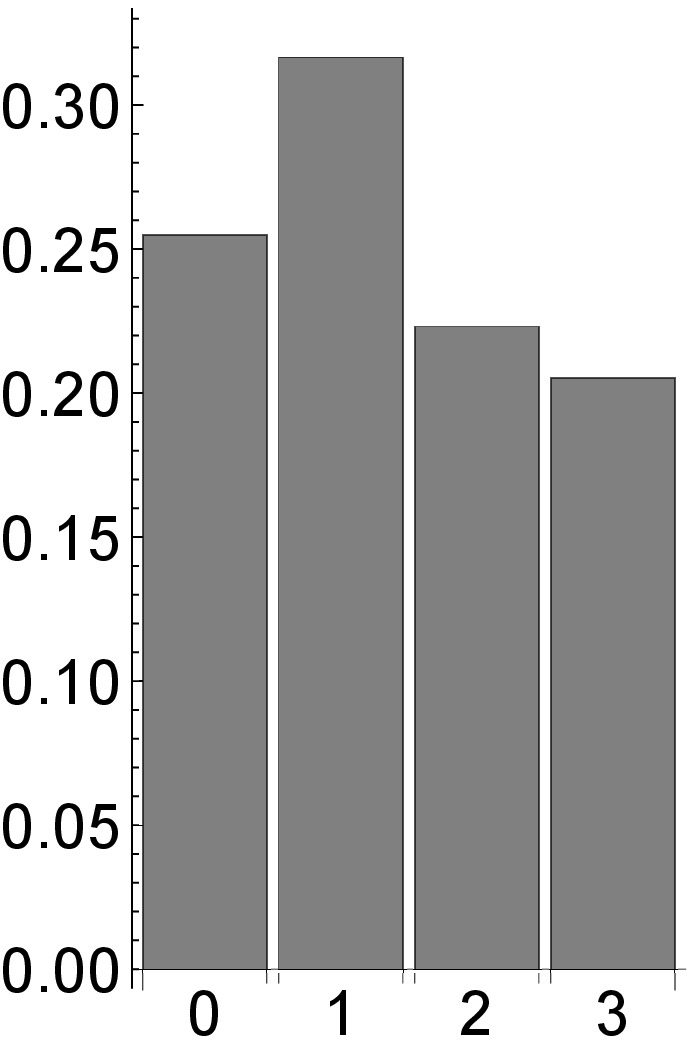}
\begin{picture}(0,0)
\put(-0.2,0.1){ $k$}
\put(-3.4,4.2){ $\P(B_2=k)$}
\end{picture}
\\
(b)&&(c)
\end{tabular}
\end{center}
\caption{The arrival pattern in Section \ref{s32}, and the probability distributions of the two platoon sizes $B_1$ and $B_2$. (a) Graphical representation of the arrival pattern; (b) Distribution of platoon 1 size; (c) Distribution of platoon 2 size.}
\label{fig:platoonsizes}
\end{figure}

In this example we consider a generalized FCTL queue with a fixed-cycle of $c=g+r=10+10=20$ time slots. Vehicles may arrive in time slots 1--10 and 16--18. The pgf of the joint distribution of $Y_1, Y_2, \dots, Y_c$ is given by
\begin{equation}
Y(z_1, z_2, \dots, z_{20})=O^{(1)}(z_1, z_2, \dots, z_{10})O^{(2)}(z_{16}, z_{17}, z_{18}),
\end{equation}
where
\begin{align*}
O^{(1)}(z_1, z_2, \dots, z_{10})
=&\ 0.0476 \ee^{0.3(z_1+z_2+z_3+z_4+z_5+z_6+z_7+z_8+z_9+z_{10}-10)}\\
& +0.107z_1\ee^{0.3(z_2+z_3+z_4+z_5+z_6+z_7+z_8+z_9+z_{10}-9)}\\
& +0.143z_1z_2\ee^{0.3(z_3+z_4+z_5+z_6+z_7+z_8+z_9+z_{10}-8)}\\
& +0.151z_1z_2z_3\ee^{0.3(z_4+z_5+z_6+z_7+z_8+z_9+z_{10}-7)}\\
&+0.138z_1z_2z_3z_4\ee^{0.3(z_5+z_6+z_7+z_8+z_9+z_{10}-6)}\\
&+0.114z_1z_2z_3z_4z_5\ee^{0.3(z_6+z_7+z_8+z_9+z_{10}-5)}\\
&+0.0887z_1z_2z_3z_4z_5z_6\ee^{0.3(z_7+z_8+z_9+z_{10}-4)}\\
&+0.0657z_1z_2z_3z_4z_5z_6z_7\ee^{0.3(z_8+z_9+z_{10}-3)}\\
&+0.0470z_1z_2z_3z_4z_5z_6z_7z_8\ee^{0.3(z_9+z_{10}-2)}\\
&+0.0328z_1z_2z_3z_4z_5z_6z_7z_8z_9\ee^{0.3(z_{10}-1)}\\
&+0.0655z_1z_2z_3z_4z_5z_6z_7z_8z_9z_{10},\\
O^{(2)}(z_1, z_2, z_{3})
=&\ 0.255 \ee^{0.075(z_1+z_2+z_3-3)}+0.317z_1\ee^{0.075(z_2+z_3-2)}\\
& +0.223z_1 z_2 \ee^{0.075(z_3-1)}+0.205z_1z_2z_3.
\end{align*}
Note that the pgfs $O^{(1)}(\cdot)$ and $O^{(2)}(\cdot)$ both have the form of the joint output pgf in \eqref{eqn:jointoutputindependent}. Each of them can be interpreted as the output of an upstream intersection with independent Poisson arrivals with intensities respectively 0.3 and 0.075 and green periods of respectively ten and three time slots. The probabilities $0.0476, 0.107, \dots, 0.0655$ and $0.255, \dots, 0.205$ correspond to the probability mass functions of the lengths of the effective green periods of the two upstream flows. Equivalently, one can also consider them as the  distributions of the platoon sizes $B_1$ and $B_2$, see Figure \ref{fig:platoonsizes}.

\paragraph{Network example without side traffic, Subsection \ref{sect:networkexamples}.}

In this network all ten intersections have synchronized settings with $c=g+r=10+10=20$. The first intersection has external arrivals from a Poisson process with intensity $\lambda=0.45$,
\[
Y^{(1)}(z_1, z_2, \dots,z_{20})=\ee^{0.45(z_1+z_2+\dots+z_{20}-20)}.
\]
Let $\rho=\lambda c/g$ denote the occupation rate, which can be interpreted as the fraction of utilized capacity at each intersection. Take $\lambda=0.45$ and hence $\rho=0.9$.

The arrival patterns at the downstream intersections depend on the delay parameter $d$. For example, when $d=5$, we have
\begin{align*}
Y^{(2)}(z_1, z_2, \dots, z_{20})
=&\ 0.0052 \ee^{0.45(z_6+z_7+z_8+z_9+z_{10}+z_{11}+z_{12}+z_{13}+z_{14}+z_{15}-10)}\\
& +0.015z_6\ee^{0.45(z_7+z_8+z_9+z_{10}+z_{11}+z_{12}+z_{13}+z_{14}+z_{15}-9)}\\
& +0.028z_6z_7\ee^{0.45(z_8+z_9+z_{10}+z_{11}+z_{12}+z_{13}+z_{14}+z_{15}-8)}\\
& +0.039z_6z_7z_8\ee^{0.45(z_9+z_{10}+z_{11}+z_{12}+z_{13}+z_{14}+z_{15}-7)}\\
&+0.048z_6z_7z_8z_9\ee^{0.45(z_{10}+z_{11}+z_{12}+z_{13}+z_{14}+z_{15}-6)}\\
&+0.054z_6z_7z_8z_9z_{10}\ee^{0.45(z_{11}+z_{12}+z_{13}+z_{14}+z_{15}-5)}\\
&+0.057z_6z_7z_8z_9z_{10}z_{11}\ee^{0.45(z_{12}+z_{13}+z_{14}+z_{15}-4)}\\
&+0.058z_6z_7z_8z_9z_{10}z_{11}z_{12}\ee^{0.45(z_{13}+z_{14}+z_{15}-3)}\\
&+0.057z_6z_7z_8z_9z_{10}z_{11}z_{12}z_{13}\ee^{0.45(z_{14}+z_{15}-2)}\\
&+0.055z_6z_7z_8z_9z_{10}z_{11}z_{12}z_{13}z_{14}\ee^{0.45(z_{15}-1)}\\
&+0.583z_6z_7z_8z_9z_{10}z_{11}z_{12}z_{13}z_{14}z_{15}.\\
\end{align*}
Again, ${0.0052,0.015,\dots,0.583}$ represent the probabilities that the platoon from the first intersection consists of, respectively, 0, 1, 2, \dots, 10 vehicles. Due to the fact that the first intersection operates close to its maximum capacity, the maximum platoon size of ten vehicles is most likely to occur.

\paragraph{Network example with side traffic, Subsection \ref{sect:networkexamples}.}

In this network we have ten queues in the main flow, $Q^{(1,0)}, \dots, Q^{(10,0)}$, with synchronized settings $c^{(i,0)}=g^{(i,0)}+r^{(i,0)}=10+10=20$ for $i=1,2,\dots,10$. The first intersection has external arrivals from a Poisson process with intensity $\lambda=0.15$,
\[
Y^{(1,0)}(z_1, z_2, \dots,z_{20})=\ee^{0.15(z_1+z_2+\dots+z_{20}-20)}.
\]
For $i=2,3, \dots, 10$, the input of $Q^{(i,0)}$ consists of the delayed output of $Q^{(i-1,0)}$ and $Q^{(i-1,1)}$, see Figure \ref{fig:networksidetraffic}.

The minor flows $Q^{(1,1)}, \dots, Q^{(9,1)}$ have cycle lengths of twenty time slots, which makes it possible to synchronize them with the major flows. However, they have different green periods of \emph{three} time slots taking place in time slots 16, 17, and 18 to avoid conflicts with the main traffic flow. The input to each of these minor flows is a Poisson process with intensity $\lambda^{(i,1)}=1/30$, for $i=1,2,\dots,9$,
\[
Y^{(i,1)}(z_1, z_2, \dots,z_{20})=\ee^{\frac{1}{30}(z_1+z_2+\dots+z_{20}-20)}.
\]
We use $\rho^{(i,0)}=\E[Y^{(i,0)}]/g^{(i,0)}$ to denote the occupation rate of the main flow of intersection $i$. Due to the inflow of side traffic, the occupation rate of the main flows increases linearly from $\rho^{(1,0)}=3/10$ to $\rho^{(10,0)}=9/10$.

These are all the required input values. The performance of the network is determined in Section \ref{sect:networkexamples} using Algorithm \ref{alg:network}.

\bibliographystyle{plain}


\end{document}